\documentclass[aap,MSNbibl,seceqn,dvips]{arximspdf}
\usepackage{mathrsfs}
\usepackage{graphicx}

\doi{10.1214/12-AAP867} 
\volume{23}
\issue{3}
\pubyear{2013}
\firstpage{1148}
\lastpage{1187}

\makeatletter

\newcommand{\rrvert}{\vert}
\newcommand{\llvert}{\vert}
\newcommand{\eqref}[1]{(\ref{#1})}
\newtheorem{theorem}{Theorem}
\newtheorem{lemma}{Lemma}
\newtheorem{proposition}{Proposition}
\newtheorem{corollary}{Corollary}
\newproclaim{remark}{Remark}
\newproclaim{remarks}{Remarks}
\newproclaim{assumptions}{Assumptions}

\def\ee{\mathrm{e}}
\def\d{\mathrm{d}}

\def\eps{\varepsilon}
\def\R{\mathbb{R}}
\def\Z{\mathbb{Z}}
\def\N{\mathbb{N}}

\newcommand{\indic}[1]{\mathbf{1}_{\{#1\}}}
\def\ind{\mathbf{1}}

\def\taut{{\widetilde{\tau}}}

\def\Var{\operatorname{Var}}

\def\bE{\mathbb{E}}
\def\bP{\mathbb{P}}

\def\bPp{\bP^{\geq0}}
\def\Pp{P^{\geq0}}
\def\bEp{\bE^{\geq0}}
\def\Ep{E^{\geq0}}

\def\nonoverlap{\operatorname{NO}(n)}
\def\overlap{\operatorname{NO}(n)^c}

\newcommand{\tautx}[1]{\taut^{(#1)}}

\newcommand{\limites}[2]{\mathop{\longrightarrow}^{#1}_{#2}}

\newcommand{\defeq}{:=}
\newcommand{\defqe}{=:}

\newcommand{\eg}{\mathbf{e}}
\newcommand{\fg}{\mathbf{f}}
\newcommand{\egb}{\bar{\mathbf{e}}}
\newcommand{\W}{W^1_\omega}
\newcommand{\WZ}{W^1_{(\omega,\widehat{Z})}}

\newcommand{\Sg}{>}
\newcommand{\Sl}{<}

\newcommand{\hP}{\mathfrak{h}} 

\makeatother

\begin{document}
\begin{frontmatter}

\title{Quenched limits for the fluctuations of transient random walks
in random environment on $\Z$\thanksref{T1}}
\thankstext{T1}{Supported by the French ANR projects MEMEMO and MEMEMO2.}
\runtitle{Quenched limits for transient RWRE}

\begin{aug}
\author[a]{\fnms{Nathana\"el}~\snm{Enriquez}\ead[label=e1]{nenriquez@u-paris10.fr}},
\author[b]{\fnms{Christophe}~\snm{Sabot}\corref{}\ead[label=e2]{sabot@math.univ-lyon1.fr}},
\author[c]{\fnms{Laurent}~\snm{Tournier}\ead[label=e3]{tournier@math.univ-paris13.fr}}
\and
\author[d]{\fnms{Olivier}~\snm{Zindy}\ead[label=e4]{olivier.zindy@upmc.fr}}

\runauthor{Enriquez, Sabot, Tournier and Zindy}
\affiliation{Universit\'e Paris 10 and Universit\'e Paris 6,
Universit\'e Lyon 1, Universit\'e~Paris~13 and Universit\'e Paris 6}
\address[a]{N. Enriquez\\
Laboratoire Modal'X\\
Universit\'e Paris 10\\
200 Avenue de la R\'epublique\\
92000 Nanterre\\
France\\
and\\
Laboratoire de Probabilit\'es\\
\quad et Mod\`eles Al\'eatoires\\
CNRS UMR 7599\\
Universit\'e Paris 6\\
4 place Jussieu\\
75252 Paris Cedex 05\\
France\\
\printead{e1}}

\address[b]{C. Sabot\\
Institut Camille Jordan\\
CNRS UMR 5208\\
Universit\'e de Lyon\\
Universit\'e Lyon 1\\
43 Boulevard du 11 novembre 1918\\
69622 Villeurbanne Cedex\\
France\\
\printead{e2}}

\address[c]{L. Tournier\\
LAGA\\
CNRS UMR 7539\\
Universit\'e Paris 13\\
99 avenue Jean-Baptiste Cl\'ement\\
93430 Villetaneuse\\
France\\
\printead{e3}}

\address[d]{O. Zindy\\
Laboratoire de Probabilit\'es\\
\quad et Mod\`eles Al\'eatoires\\
CNRS UMR 7599\\
Universit\'e Paris 6\\
4 place Jussieu\\
75252 Paris Cedex 05\\
France\\
\printead{e4}\hspace*{34pt}}

\end{aug}

\received{\smonth{5} \syear{2011}}
\revised{\smonth{3} \syear{2012}}

%
\begin{abstract}
We consider transient nearest-neighbor random walks in random
environment on $\Z.$
For a set of environments whose probability is converging to~1 as time
goes to infinity, we describe the fluctuations of the hitting time of a
level $n$, around its mean, in terms of an explicit function of the
environment. Moreover, their limiting law is described using a Poisson
point process whose intensity is computed. This result can be
considered as the quenched analog of the classical result of Kesten,
Kozlov and Spitzer [\textit{Compositio Math}. \textbf{30} (1975) 145--168].
\end{abstract}

%


\end{frontmatter}

\section{Introduction}
\label{sintro}

Random walks in a one-dimensional random environment were first
introduced in the late sixties as a toy model for DNA replication. The
recent development of micromanipulation techniques such as DNA
unzipping has raised a renewed interest in this model in genetics and
biophysics; cf., for instance,~\cite{cocco} where it is involved in a
DNA sequencing procedure. Its mathematical study was initiated by
Solomon's 1975 article~\cite{solomon} characterizing the transient and
recurrent regimes and proving a strong law of large numbers. A salient
feature emerging from this work was the existence of an intermediary
regime where the walk is transient with a zero asymptotic speed, in
contrast with the case of simple random walks. Shortly after, Kesten,
Kozlov and Spitzer~\cite{kks} precised this result by giving limit laws
in the transient regime. When suitably normalized, the (properly
centered) hitting time of site $n$ by the random walk was proved to
converge toward a stable law as $n$ tends to infinity, which implies a
limit law for the random walk itself. In particular, this entailed that
the ballistic case (i.e., with positive speed) further decomposes into
a diffusive and a subdiffusive regime.

Note that these results, except when they deal with almost sure
statements, concern only the annealed behavior. When dealing with
applications, what we call the medium is usually fixed during the
experiment (e.g., the DNA sequence), and we are naturally led to
consider the quenched behavior of the walk. The first results in this
direction by Peterson and Zeitouni~\cite{peterson-zeitouni} and
Peterson~\cite{peterson} were unfortunately negative results, saying
that, for almost all environment, the laws of the fluctuations of the
walk along the time have several accumulation points. However, it was
shown by three of the authors in~\cite{aging}, that, in the case of
transient walks having 0 asymptotic speed, one can get some quenched
localization result by slightly relaxing the point of view. Namely, for
a set of media whose probability converges to 1 as time goes to
infinity, the law of the (suitably normalized) position of the walk is
getting close to a discrete probability measure whose weights and
support are expressed in terms of the environment. In the same spirit,
we focus in this work on the quenched fluctuations of hitting times in
the case of a general transient subdiffusive random walk in random environment.

Adopting Sinai's now famous description of the medium by a potential
\cite{sinai}, we introduce a notion of valley. We then prove that the
fluctuations of the hitting time of $x$ around its expectation mainly
come from the times spent crossing a very small number of deep
potential wells. Since these wells are well apart, their crossing times
are almost independent. Moreover, it is shown that the laws of these
crossing times are well approximated by exponential variables whose
expectations are functions of the environment, functions which in turn
happen to be closely related to the classical Kesten renewal series.

Thus, our main result states that the law of the difference of a
hitting time with its expectation is close to the law of a sum of
centered exponential variables which are weighted by heavy-tailed
functions of the environment. This makes it possible to describe their
law in terms of a Poisson point process whose intensity is explicitly computed.

To make the exposition clearer, we first present the main results and
notation (Section~\ref{shyp+thm}) and defer to Section \ref
{secsketch} the more precise description of the organization of the
paper along with a sketch of the proof.

\section{Notation and main results}
\label{shyp+thm}

Let $\omega\defeq(\omega_x, x \in\Z)$ be a family of i.i.d.
random variables taking values in $(0,1)$, which stands for the random
environment. Let $\Omega\defeq(0,1)^\Z$ and denote by $P$ the
distribution of $\omega$ (on $\Omega$) and by $E$ the corresponding
expectation. Conditioning on $\omega$ (i.e., choosing an environment),
we define the random walk in random environment $X\defeq(X_t, t\in
\N
)$ starting from $x\in\Z$ as a nearest-neighbor random walk on $\Z$
with transition\vadjust{\goodbreak} probabilities given by $\omega$: if we denote by
$P_{x,\omega}$ the law of the Markov chain $(X_t, t\ge0)$ defined by
$P_{x,\omega} ( X_{0} = x ) =1$ and
\[
P_{x,\omega} ( X_{t+1} = z | X_t =y ) \defeq
\cases{
\omega_y, & \quad$\mbox{if } z=y+1,$
\vspace*{2pt}\cr
1-\omega_y, &\quad$\mbox{if } z=y-1,$
\vspace*{2pt}\cr
0, &\quad$\mbox{otherwise,} $}
\]
then the joint law of $(\omega,X)$ is $\bP_x(\d\omega,\d X)\defeq
P_{x,\omega}(\d X)P(\d\omega)$.
For convenience, we let $\bP\defeq\bP_0$. We refer to~\cite{zeitouni}
for an overview of results on random walks in random environment.
An important role is played by the sequence of variables
%
%
\begin{equation}
\label{defrho}
\rho_x\defeq\frac{1-\omega_x}{\omega_x},\qquad x \in\Z.
\end{equation}
We will make the following assumptions in the rest of this paper.
%
%
\begin{assumptions*}
\begin{longlist}[(a)]
\item[(a)] There exists $0<\kappa<2$ for which $E
[\rho_0^{\kappa} ]=1$ and $E [ \rho_0^{\kappa} \log^+ \rho_0
]<\infty$;
\item[(b)] The distribution of $\log\rho_0$ is nonlattice.
\end{longlist}
\end{assumptions*}
Let us recall here that, under assumptions (a) and (b), Kesten,
Kozlov and Spitzer~\cite{kks} proved a limit theorem toward a stable
law of index $\kappa,$ whose scaling parameter is obtained in \cite
{limitlaws} for the sub-ballistic case and in~\cite{stable} for the
ballisitic case.

We now introduce the hitting time $\tau(x)$ of site $x$ for the
random walk $(X_t, t \ge0),$
\[
\label{hittimerw} \tau(x)\defeq\inf\{ t \ge0\dvtx X_t=x \},\qquad
x \in
\Z,
\]
and the inter-arrival time $\tau(x,y)$ between sites $x$ and $y$ by
\[
\label{hittimerwia} \tau(x,y) \defeq\inf\{ t \ge0\dvtx X_{\tau
(x)+t}=y \},\qquad
x,y \in\Z.
\]

Following Sinai~\cite{sinai} (in the recurrent case), and more recently
the study of the case $0<\kappa<1$ in~\cite{limitlaws}, we define a
notion of potential that enables us to visualize where the random walk
spends most of its time.

The potential, denoted by $V= (V(x), x\in\Z)$, is a function of the
environment~$\omega$ defined by $V(0)=0$ and $\rho_x=\ee
^{V(x)-V(x-1)}$ for every
$x\in\Z$, that~is,
\[
V(x) \defeq\cases{
\displaystyle\sum
_{1\leq y\leq x} \log\rho_y, &\quad$\mbox{if } x \ge1,$
\vspace*{2pt}\cr
0, &\quad$\mbox{if } x=0,$
\vspace*{2pt}\cr
\displaystyle-\sum_{x<y\leq0}^0 \log\rho_y,
&\quad$\mbox{if } x\le-1,$}
\]
where the $\rho_y$'s are defined in \eqref{defrho}. Under hypothesis
(a), Jensen's inequality gives $E[\log\rho_0^\kappa]\leq\log
E[\rho_0^\kappa]=0$, and hypothesis (b) excludes the equality
case $\rho_0=1$ a.s., hence, $E[\log\rho_0]< 0$ and thus $V(x)\to
\mp\infty$ a.s. when $x\to\pm\infty$.\vadjust{\goodbreak}

The potential is subdivided into pieces, called ``excursions,'' by its
weak descending ladder epochs $(e_p)_{p\geq0}$ defined by $e_0\defeq
0$ and
%
%
\begin{equation}
\label{eqndefei} e_{p+1} \defeq\inf\bigl\{ x > e_p\dvtx
V(x) \le V(e_p)\bigr\},\qquad p \ge0.
\end{equation}
The number of excursions before $x>0$ is
%
%
\begin{equation}
\label{eqndefnx} n(x)\defeq\max\{p\dvtx e_p\leq x\}.
\end{equation}

Moreover, let us introduce the constant $C_K$ describing the tail of
Kesten's renewal series $R\defeq\sum_{x\geq0}\rho_0\cdots\rho
_x=\sum_{x\geq0}\ee^{V(x)}$ (see~\cite{kesten73}) that plays a
crucial role
in this work:
\[
P(R>t)\sim C_K t^{-\kappa},\qquad t \to\infty.
\]
Note that at least two probabilistic representations are available to
compute $C_K$ numerically, which are equally efficient. The first one
was obtained by Goldie~\cite{goldie} and a second one was obtained in
\cite{renewal}.

Finally, recall the definition of the Wasserstein metric $W^1$ between
probability measures $\mu,\nu$ on $\R$:
\[
W^1(\mu,\nu) \defeq\mathop{\inf_{(X,Y):}}_{X\sim\mu, Y\sim\nu
}
E\bigl[|X-Y|\bigr],
\]
where the infimum is taken over all couplings $(X,Y)$ with marginals
$\mu$ and $\nu$. We will denote by $\W(X,Y)$ the $W^1$ distance between
the laws of random variables~$X$ and $Y$ conditional to $\omega$, that
is, between the ``quenched distributions'' of $X$ and $Y$.

Let us emphasize that the following results, which describe the
quenched law of $\tau(x)$ in terms of the environment, can be stated in
different ways, depending on the applications we have in mind, either
practical or theoretical. We give two variants and mention that the
following results hold for any $\kappa\in(0,2)$ (so that the
sub-ballistic regime is also included, even though a finer study was
led for $\kappa\in(0,1)$ in~\cite{aging}).

%
\begin{theorem} \label{tquenched2}
Under assumptions \textup{(a)} and \textup{(b)} we have
\[
\W\Biggl(\frac{\tau(x)-E_\omega[\tau(x)]}{x^{1/\kappa}},\frac
{1}{x^{1/\kappa}}\sum_{p=0}^{n(x)-1}
E_\omega\bigl[\tau(e_p,e_{p+1})\bigr]
\egb_p \Biggr)\limites{P\mathrm{\mbox{-}probability}} {x}0,
\]
with $\egb_p:=\eg_p-1$, where $(\eg_p)_p$ are i.i.d. exponential
random variables of parameter~$1$ independent of $\omega$; the terms
$E_\omega[\tau(e_p,e_{p+1})]$ can be made explicit $[$see $(\ref
{eqnzeitounie})$ in the Preliminaries$]$, and $n(x)$ may be replaced
by $ \lfloor\frac{x}{E[e_1]} \rfloor$.
\end{theorem}

%
\begin{theorem} \label{tquenched3}
Under assumptions \textup{(a)} and \textup{(b)}, for every $\delta>0$
and $\eps>0$, if~$x$ is large enough, we may enlarge the probability
space so as to introduce i.i.d. random variables $\widehat
{Z}=(\widehat{Z}_{p})_{p\geq0}$ such that
%
%
\begin{equation}
\label{eqnqueuezh} P(\widehat{Z}_p>t)\sim2^\kappa
C_U t^{-\kappa},\qquad t\to\infty,\vadjust{\goodbreak}
\end{equation}
where $C_U\defeq E[\rho_0^\kappa\log\rho_0]E[e_1](C_K)^2$, and
\[
P \Biggl(\WZ\Biggl(\frac{\tau(x)-E_\omega[\tau(x)]}{x^{1/\kappa
}},\frac
{1}{x^{1/\kappa}}\sum
_{p=1}^{\lfloor x/{E[e_1]}\rfloor} \widehat{Z}_p
\egb_p \Biggr)>\delta\Biggr)<\eps,
\]
with $\egb_p:=\eg_p-1$, where $(\eg_p)_p$ are i.i.d. exponential
random variables of parameter~$1$ independent of $\widehat{Z}$, and
$\WZ
(X,Y)$ denotes the $W^1$ distance between the law of $X$ given $\omega$
and the law of $Y$ given $\widehat{Z}$. 
\end{theorem}


By a classical result (cf.~\cite{durrett}, page 152, or~\cite{resnick},
page 138, for a general statement), the set $\{n^{-1/\kappa}\widehat
{Z}_p | 1\leq p\leq n\}$ converges toward a Poisson point process
of intensity $2^\kappa C_U \kappa x^{-(\kappa+1)}\,\d x$. It is
therefore natural to expect the following corollary.

%
\begin{corollary}\label{coroll}
Under assumptions \textup{(a)} and \textup{(b)} we have
\[
\mathscr{L}\biggl(\frac{\tau(x)-E_\omega[\tau
(x)]}{x^{1/\kappa}} \Big\vert\omega\biggr)
\limites{W^1} {x} \mathscr{L} \Biggl(\sum
_{p=1}^\infty\xi_p \egb_p
\Big\rrvert(\xi_p)_{p\geq1} \Biggr) \qquad\mbox{in law},
\]
where the convergence is the convergence in law on the $W^1$ metric
space of probability measures on $\R$ with finite first moment, and
$(\xi_p)_{p\geq1}$ is a Poisson point process of intensity $\lambda
\kappa u^{-(\kappa+1)}\,\d u$ where
\[
\lambda\defeq\frac{2^\kappa C_U}{E[e_1]}=2^\kappa\kappa E\bigl[
\rho_0^\kappa\log\rho_0\bigr]C_K^2,
\]
$\egb_p:=\eg_p-1$ where $(\eg_p)_p$ are i.i.d. exponential random
variables of parameter~$1$, and the two families are independent of
each other. In the case $\kappa=1$, $\lambda=\frac{2}{E[\rho_0\log
\rho_0]}$, and in the case where $\omega_0$ has a distribution
$\operatorname{
Beta}(\alpha,\beta)$, with $0\Sl\alpha-\beta\Sl2$,
\[
\lambda=2^{\alpha-\beta}\frac{\Psi(\alpha)-\Psi(\beta)}{(\alpha
-\beta)
B(\alpha-\beta,\beta)^2},
\]
where $\Psi$ denotes the classical digamma function
$\Psi(z)\defeq(\log\Gamma)'(z)=\frac{\Gamma'(z)}{\Gamma(z)}$ and
$B(\alpha,\beta)\defeq\int_{0}^1 x^{\alpha-1}(1-x)^{\beta-1} \,\d
x=\frac
{\Gamma(\alpha)\Gamma(\beta)}{\Gamma(\alpha+\beta)}$.
\end{corollary}

%
\begin{remarks*}
\begin{longlist}[(1)]
\item[(1)] Since the topology of convergence in $W^1$-distance is
finer than the topology of weak convergence restricted to probability
measures with finite first moment, we may replace $W^1$ by the topology
of the convergence in law in the above limit.
\item[(2)] For every $\eps>0$, the mass of $(\eps,+\infty)$ for the
measure $\mu=\lambda\kappa\frac{\d u}{u^{\kappa+1}}$ is finite so that
it makes sense to consider a decreasing ordering\vadjust{\goodbreak} $(\xi^{(k)})_{k\geq
1}$ of the Poisson process of intensity $\mu$. A change of variable
then shows that
%
%
\begin{equation}
\label{eqndescrxi} \xi^{(p)}=\lambda^{1/\kappa}(\fg_1+
\cdots+\fg_p)^{-1/\kappa
},\qquad p\geq1,
\end{equation}
$(\fg_p)_p$ being i.i.d. exponential random variables of parameter $1$.
In particular, by the law of large numbers,
%
%
\begin{equation}
\label{eqnequivxi} \xi^{(p)}\sim\lambda^{1/\kappa} p^{-1/\kappa
},\qquad
p\to\infty, \mbox{ a.s.},
\end{equation}
hence, $\sum_p(\xi_p)^2=\sum_p(\xi^{(p)})^2<\infty$ a.s. Thus, the
random series $\sum_p \xi_p \egb_p$ converges a.s. Furthermore, since
its characteristic function is also an absolutely convergent product,
its law does not depend on the ordering of the points.
\end{longlist}
\end{remarks*}

Corollary~\ref{coroll} can be easily deduced from the previous theorems.
We give a short proof of this result in Section~\ref{secproofcoro}.

While finishing writing the present article, we learned about the
article~\cite{peterson-samorodnitsky} by Peterson and Samorodnitsky
giving a result close to Corollary~\ref{coroll}. Another article \cite
{dolgopyat-goldsheid} by Dolgopyat and Goldsheid was also submitted,
that establishes a similar result (under the ellipticity condition).
Our statement, however, gives the convergence in $W^1$ instead of the
weak convergence and especially specifies the \textit{value of the
constant} $\lambda$ that appears in the intensity of the limiting
Poisson point process. Furthermore, the three proofs are rather different.

In the following, the constant $C$ stands for a positive constant large
enough, whose value can change from line to line.

\section{Sketch of the proof}
\label{secsketch}



Along the sequence $(e_p)_{p\geq0}$, hitting times decompose into
crossing times of a linear number of excursions,
\[
\tau(x)=\sum_{0\leq p<n(x)}\tau(e_p,e_{p+1})
+ \tau(e_{n(x)},x).
\]
Although these terms are very correlated, the core of the proof
consists of the fact that, as far as fluctuations are concerned, the
main contribution only comes from a logarithmic subfamily of
asymptotically i.i.d. terms which correspond to so-called ``high
excursions'' (or ``deep valleys''). This property (stemming from the
fact that the random variables $E_\omega[\tau(e_p,e_{p+1})]$ are
heavy-tailed) enables the proof to be divided into two parts detailed below.

\subsection{\texorpdfstring{Exit time from a deep valley (Section \protect\ref{secexitsinglevalley})}{Exit time from a deep valley (Section 5)}}
\label{subsecexitonevalley}

\begin{figure}

\includegraphics{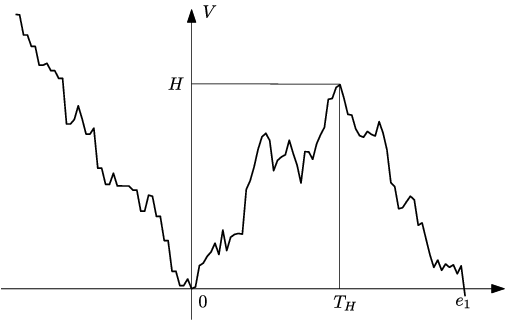}

\caption{Height of an excursion.}\label{figexcursion}
\end{figure}

The crossing time of the excursion $[e_p,e_{p+1}]$ will mainly depend
on its height
\[
H_p\defeq\max_{e_p\leq x<e_{p+1}} \bigl(V(x)-V(e_p)
\bigr).
\]
%
As $p$ grows, the law of the potential $V$ viewed from $e_p$ converges to
$\Pp\defeq P( \cdot| \forall x\leq0, V(x)\geq0)$, and therefore
the time $\tau(e_p,e_{p+1})$ converges in law to $\tau(e_1)$ under
$\Pp
$ which we have now to study.
A classical Markov chain computation gives (cf. Section~\ref{subsecquench})
\[
E_\omega\bigl[\tau(e_1)\bigr]=\sum
_{0\leq y<e_1}\sum_{x\leq y}(2-\indic{x=y})
\ee^{V(y)-V(x)}.
\]
When $H\defeq H_0$ is large, factorizing by the largest term $2\ee^{H}$
leads to
\[
E_\omega\bigl[\tau(e_1)\bigr]\simeq2 \ee^H\sum
_x \ee^{-V(x)} \sum
_y \ee^{-(H-V(y))},
\]
where in the sums the significant terms are those indexed by values $x$
close to $0$ and values $y$ close to $T_H$; cf. Figure \ref
{figexcursion}. In particular, we have
\[
E_\omega\bigl[\tau(e_1)\bigr]\simeq2 \ee^H
M_1 M_2,
\]
where $M_1$, $M_2$ are defined by \eqref{eqdefM1M2}. Due to the
``locality'' of $M_1$ and $M_2$, a key fact from~\cite{renewal} is
that, when $H$ is large, $M_1$, $M_2$ and $H$ are asymptotically
independent and $M_1$, $M_2$ have the same law. Now, Iglehart's tail
estimate on $\ee^{H}$ [see \eqref{iglehartthm}] yields
\[
\Pp\bigl(E_\omega\bigl[\tau(e_1)\bigr]\geq t \bigr)
\sim2^\kappa C_I E\bigl[M^{\kappa}
\bigr]^2 t^{-\kappa},\qquad t\to\infty,
\]
where $C_I$ is given by \eqref{cstIglehart}. This is an important
result of~\cite{renewal}, rephrased in Lemma~\ref{lemtailz}.

%

To complete the description of the law of the crossing time of a ``high
excursion,'' we furthermore prove in Section~\ref{secexitsinglevalley}
that, for large $H$, the law of $\tau(e_1)$, given~$\omega$, is close
to an exponential law with mean $E_\omega[\tau(e_1)]$. This follows
from the fact that the number of returns to 0 before reaching $e_1$
follows a geometric law.

%

\subsection{\texorpdfstring{Deep and shallow valleys (Sections \protect\ref{seciidvalleys} and \protect\ref{secinterarrival})}
{Deep and shallow valleys (Sections 6 and 7)}}
\begin{figure}

\includegraphics{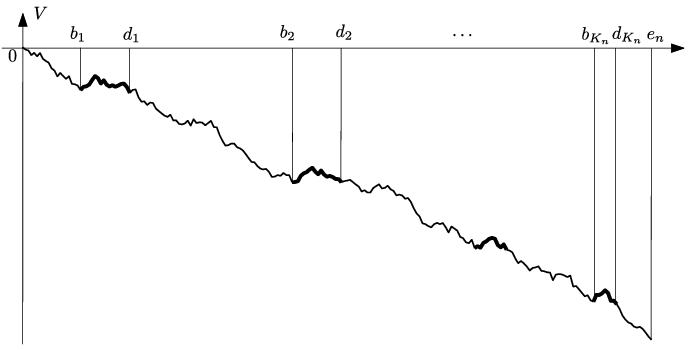}

\caption{High excursions $($in bold$)$ among the $n$ first
excursions.}\label{figexcursions}
\end{figure}

As mentioned at the beginning of the section, we try to focus the study
on the crossing times of high excursions. To this aim, we introduce a
critical height $h_n$, adapted to the space scale $n$, defined by
\[
\label{hcritic} h_n\defeq\frac{1}{\kappa}\log n - \log\log n.
\]
Then, let $(\sigma(i))_{i \ge1}$ be the sequence of the indices of the
successive excursions whose heights are greater than $h_n.$ More precisely,
\begin{eqnarray*}
\sigma(1)&\defeq&\inf\{ p \ge0 \dvtx H_p \ge h_n\},
\\
\sigma(i+1)&\defeq&\inf\bigl\{ p >\sigma(i)\dvtx H_p \ge
h_n \bigr\}, \qquad i \ge1.
\end{eqnarray*}
The \textit{high excursions} (see Figure~\ref{figexcursions}) are defined
as the restriction of the potential to $[b_i, d_i]$, where
\[
b_i\defeq e_{\sigma(i)}, \qquad{d}_i\defeq
e_{\sigma(i)+1}.
\]
Note that by Iglehart's estimate, the probability $P(H \geq h_n )$ is
asymptotically equal to $C_I \ee^{- \kappa h_n}$, hence, the number
of high excursions among the $n$ first ones,
\[
K_n\defeq\# \{ 0\leq i\leq n-1 \dvtx H_i\geq
h_n \},
\]
is of order $(\log n)^\kappa$.

%

It turns out that the crossing time $\tau(b_i,d_i)$ involves mainly the
environment between $a_i$ and $d_i$ where $a_i$ is defined
as
\[
a_i\defeq e_{\sigma(i)-D_n}
\]
and $D_n$ is chosen in such a way that $V(a_i)-V(b_i)$ is slightly
greater than $h_n$, that is,
%
%
\begin{equation}
\label{eqndefdn} D_n\defeq\biggl\lceil\frac{1+\gamma}{A \kappa
}\log n
\biggr\rceil,
\end{equation}
$\gamma>0$ being arbitrary and $A$ being equal to $E[-V(e_1)]$ if this
expectation is finite, and otherwise being an arbitrary positive real number.
The \textit{deep valleys} are defined as the restriction of the potential
to $[a_i, d_i]$.

We successively prove that:
\begin{longlist}[(1)]
\item[(1)] deep valleys are asymptotically disjoint and their exit
times $\tau(b_i,d_i)$ are asymptotically i.i.d. (Section \ref
{seciidvalleys});
\item[(2)] the contribution to fluctuations of the crossing times of
low excursions is negligible (Section~\ref{secinterarrival}).
\end{longlist}
This second point constitutes a novelty with respect to previous works
in that the contribution of the crossing times of the numerous small
excursions is not negligible with respect to $\tau(x)$ (for $1\le
\kappa<2$) but only their fluctuations are, and for this reason we
have to control their covariances.

The behavior summarized above is emphasized in the following
formulation which lies at the core of the proof: under assumptions
(a) and (b),
%
%
\begin{equation}
\label{tequenched} \W\Biggl(\frac{\tau(x)-E_\omega[\tau
(x)]}{x^{1/\kappa}},\frac
{1}{x^{1/\kappa}}\sum
_{i=1}^{K(x)} E_\omega\bigl[
\tau(b_i,d_i)\bigr]\egb_i \Biggr)
\mathop{\longrightarrow}^{P\mathrm{\mbox{-}probability}}_{x}
0,
\end{equation}
with $\egb_i:=\eg_i-1$ where $(\eg_i)_i$ are i.i.d. exponential random
variables of parameter~1, independent of $\omega$, and $K(x):=K_{n(x)}$
where $n(x)$ is defined by~\eqref{eqndefnx}.

Note that the terms $E_\omega[\tau(b_i,d_i)]$ can be made explicit [see
\eqref{eqnzeitounie}]. Hence, this formula is well suited to derive
practical information about $\tau(x)$ which, for instance, appears as
an unzipping time in~\cite{cocco}.

\section{Preliminaries} \label{secpreliminaries}

This section is divided into three independent parts. The first part
quickly recalls a stationarity property of the potential when suitably
conditioned on $\Z_-$, which is used throughout the paper. The second
one recalls usual formulas about random walks in a one-dimensional
potential. Finally, the last part adapts the main results from \cite
{renewal} in the present context.

\subsection{Environment on the left of $0$}\label{subsecleftof0}
It will be convenient to extend the sequence $(e_p)_{p\geq0}$ to
negative indices by letting
%
%
\begin{equation}
\label{eqndefenegatifs} e_{p-1}\defeq\sup\bigl\{x< e_p
\dvtx\forall y<x, V(y)\geq V(x)\bigr\},\qquad p\leq0.
\end{equation}
The structure of the sequence $(e_p)_{p\in\Z}$ will be better
understood after Lemma~\ref{lemenegatifs}.

We accordingly extend the sequence $(H_p)_{p\geq0}$ of heights
\[
H_p \defeq\max_{e_p \le x \le e_{p+1}} \bigl(V(x)-V(e_p)
\bigr),\qquad p\in\Z.\vadjust{\goodbreak}
\]
Note that the excursions $(V(e_p+x)-V(e_p))_{0\leq x< e_{p+1}-e_p}$,
$p\geq0$, are i.i.d. Also, the intervals $(e_p,e_{p+1}], {p \in\Z},$
stand for the excursions of the potential above its past minimum,
provided $V(x)\geq0$ when $x\leq0$.

By definition, the distribution of the environment is translation
invariant. However, the distribution of the ``environment seen from
$e_p$,'' that is, of $(\omega_{e_p+x})_{x\in\Z}$, depends on $p\in
\Z$.
When suitably conditioning the environment on $\Z_-$, this problem vanishes.

Let us define the conditioned probabilities
\[
\Pp\defeq P\bigl( \cdot| \forall x\leq0, V(x)\geq0\bigr) \quad
\mbox{and}\quad\bPp
\defeq P_\omega\times\Pp(d\omega).
\]
Then the definition of $e_p$ for $p<0$ classically implies the
following useful property.

%
\begin{lemma}\label{lemenegatifs}
Under $\Pp$, the sequence $(V(e_p+x)-V(e_p))_{0\leq x\leq
e_{p+1}-e_p}$, $p\in\Z$, of excursions is i.i.d. In particular,
the sequence $(V(e_p+x)-V(e_p))_{x\in\Z}$ of potentials $[$and thus the
sequence $(\omega_{e_p+x})_{x\in\Z}$, $p\in\Z$, of environments$]$ is
stationary under $\Pp$.
\end{lemma}

\subsection{Quenched formulas}\label{subsecquench}

We recall here a few Markov chain formulas that are of repeated use in
the paper.

\textit{Quenched exit probabilities}.
For any $a\leq x\leq b$ (see~\cite{zeitouni}, formula (2.1.4))
\begin{equation}\label{eqnzeitounip}
P_{x,\omega}\bigl(\tau(b)<\tau(a)\bigr)=\frac{\sum_{a\leq y<
x}\ee
^{V(y)}}{\sum_{a\leq y<b}\ee^{V(y)}}.
\end{equation}
In particular,
\begin{equation}\label{eqnzeitounip3}
P_{a+1,\omega}\bigl(\tau(a)=\infty\bigr)= \biggl(\sum_{y\geq
a}\ee^{V(y)-V(a)} \biggr)^{-1}.
\end{equation}
Thus, $P_{0,\omega}(\tau(1)=\infty)= (\sum_{x\le0}\ee^{V(x)}
)^{-1} =0$, $P$-a.s. because $V(x)\to+\infty$ a.s. when $x\to
-\infty
$, and $P_{1,\omega}(\tau(0)=\infty)= (\sum_{x\ge0}\ee^{V(x)}
)^{-1} >0$, $P$-a.s. by the root test (using $E[\log\rho_0]<0$). This
means that $X$ is transient to $+\infty,$ $\bP$-a.s.

\textit{Quenched expectation}. For any $a<b$, $P$-a.s.
(cf.~\cite{zeitouni})
%
%
\begin{equation}
E_{a,\omega}\bigl[\tau(b)\bigr] = \sum_{a\leq y<b}
\sum_{x\leq y}\alpha_{xy}\ee^{V(y)-V(x)},
\label{eqnzeitounie}
\end{equation}
where $\alpha_{xy}=2$ if $x<y$, and $\alpha_{yy}=1$. Thus, we have
\begin{equation}\label{eqnzeitounie2}
E_{a,\omega}\bigl[\tau(b)\bigr]
\leq2\sum_{a\leq y<b}\sum_{x\leq y}\ee^{V(y)-V(x)}
\end{equation}
and, in particular,
\begin{equation}\label{eqnzeitounie3}
E_{a,\omega}\bigl[\tau(a+1)\bigr]=1+2\sum_{x< a}\ee
^{V(a)-V(x)}\leq2\sum
_{x\leq
a}\ee^{V(a)-V(x)}.
\end{equation}
\textit{Quenched variance}. For any $a<b$, $P$-a.s. (cf.
\cite{alili} or~\cite{goldsheid})
%
%
\begin{eqnarray}
\Var_{a,\omega}\bigl(\tau(b)\bigr) & =& 4\sum_{a\leq y<b}
\sum_{x\leq y}\ee^{V(y)-V(x)}\bigl(1+
\ee^{V(x-1)-V(x)}\bigr)
\nonumber
\\[-8pt]
\\[-8pt]
\nonumber
&&\hspace*{49pt}{}\times \biggl(\sum_{z<x}
\ee^{V(x)-V(z)} \biggr)^2, \label{eqnalili}
\end{eqnarray}
from where we get, after expansion, change of indices and addition of a
few terms,
%
%
\begin{equation}
\Var_{a,\omega}\bigl(\tau(b)\bigr) \leq16\sum
_{a\leq y< b}\sum_{z'\leq z\leq x\leq y}
\ee^{V(y)+V(x)-V(z)-V(z')}. \label{eqnalili2}
\end{equation}

\subsection{Renewal estimates}
\label{subsecrenewal}

In this section we recall and adapt results from~\cite{renewal}, which
are very useful to finely bound the expectations of exponential
functionals of the potential.

Let us first observe that hypothesis (a) implies that $e_1$ is
exponentially integrable. Indeed, for all $x\in\N$, for any $\lambda
>0$, $P(e_1>x)\leq P(V(x)>0)=P(\ee^{\lambda V(x)}>1)\leq E[\ee
^{\lambda
V(x)}]= E[\rho_0^\lambda]^x$, and $E[\rho_0^\lambda]<1$ for any
$0<\lambda<\kappa$ by convexity of $s\mapsto E[\rho_0^s]$.

Let $R_-\defeq\sum_{x\leq0}\ee^{-V(x)}$. Then, Lemma 3.2 from \cite
{renewal} proves that
%
%
\begin{equation}
\label{eqnr-} \Ep[R_-]<\infty,
\end{equation}
and that more generally all the moments of $R_-$ are finite under $\Pp$.

The study of ``high excursions'' involves the following key result of
Iglehart~\cite{igle} which gives the tail probability of $H$ (recall
$H\defeq H_0$), namely,
%
%
\begin{equation}
\label{iglehartthm} P (H \geq h) \sim C_I \ee^{- \kappa h},\qquad h
\to
\infty,
\end{equation}
where
%
%
\begin{equation}
\label{cstIglehart} C_I\defeq{\frac{(1-E[\ee^{\kappa
V(e_1)}])^2}{\kappa E[\rho_0^\kappa\log
\rho_0]E[e_1]}}.
\end{equation}

Let us define
\[
T_H\defeq\min\bigl\{x\geq0\dvtx V(x)= H\bigr\}
\]
and
%
%
\begin{equation}
\label{eqdefM1M2} M_1\defeq\sum_{x<T_H}
\ee^{-V(x)},\qquad M_2\defeq\sum_{0\leq
x<e_1}
\ee^{V(x)-H}.
\end{equation}

Let $Z\defeq M_1 M_2 \ee^H$. Theorem 2.2 (together with Remark A.1) of
\cite{renewal} proves that
%
%
\begin{equation}
\label{eqnqueuezcondition} \Pp(Z>t, H=S)\sim C_U t^{-\kappa},\qquad
t\to\infty,
\end{equation}
where $C_U$ was defined after \eqref{eqnqueuezh}; cf. also the sketch
in Section~\ref{subsecrenewal} for heuristics.
While the condition $\{H=S\}$ was natural in the context of \cite
{renewal}, we will need to remark that we may actually drop it.

%
\begin{lemma}\label{lemtailz}
We have
\[
\Pp(Z>t)\sim C_U t^{-\kappa},\qquad t\to\infty.
\]
\end{lemma}

The proof of this lemma is postponed to Appendix~\ref
{subsecprooflemma-renewalestimate}.
We will often need moments involving
\[
M'_1\defeq\sum_{x<e_1}
\ee^{-V(x)},
\]
instead of $M_1(\leq M'_1)$. The next result is an adaptation of Lemma
4.1 from~\cite{renewal} to the present situation, together with \eqref
{iglehartthm}, with a novelty coming from the difference between $M'_1$
and $M_1$.

%
\begin{lemma}\label{lemrenewalestimates}
For any $\alpha,\beta,\gamma\geq0$, there is a constant $C$ such that,
for large $h>0$,
%
%
\begin{equation}
\label{eqnm1m2inf} \Ep\bigl[\bigl(M'_1
\bigr)^\alpha(M_2)^\beta\ee^{\gamma H}|H< h
\bigr]\leq\cases{
C, & \quad$\mbox{if }\gamma<\kappa,$
\vspace*{2pt}\cr
C h,&\quad$\mbox{if }\gamma=\kappa,$
\vspace*{2pt}\cr
C \ee^{(\gamma-\kappa)h},&\quad$\mbox{if }\gamma>\kappa,$}
\end{equation}
and, if $\gamma<\kappa$,
%
%
\begin{equation}
\label{eqnm1m2sup} \Ep\bigl[\bigl(M'_1
\bigr)^\alpha(M_2)^\beta\ee^{\gamma H}|H\geq h
\bigr]\leq C\ee^{\gamma h}.
\end{equation}
\end{lemma}

The proof of this lemma is technical and therefore postponed to
Appendix~\ref{subsecprooflemma-renewalestimate}.
Let us now give an important application of Lemma~\ref{lemrenewalestimates}.

%
\begin{lemma} \label{lemboundEe}
We have, for all $h>0$, if $0<\kappa<1$,
%
%
\begin{equation}
\label{eqnboundEe} \Ep\bigl[E_\omega\bigl[\tau(e_1)\bigr]
\indic{H<h} \bigr]\leq C\ee^{(1-\kappa)h}
\end{equation}
and, if $0<\kappa<2$,
%
%
\begin{equation}
\label{eqnboundEe2} \Ep\bigl[E_\omega\bigl[\tau(e_1)
\bigr]^2\indic{H<h} \bigr]\leq C\ee^{(2-\kappa)h}.
\end{equation}
\end{lemma}

\begin{pf}
Since, by \eqref{eqnzeitounie2}, we have $E_\omega[\tau(e_1)]\leq
2M_1'M_2\ee^H$, the result follows directly from Lemma \ref
{lemrenewalestimates}.
\end{pf}

\section{Exit time from a deep valley}
\label{secexitsinglevalley}

This section aims at proving that the quenched law of the crossing time
\[
\tau\defeq\tau(e_1)
\]
of an excursion is close to that of $E_\omega[\tau]\eg$, where $\eg
$ is
an exponential random variable independent of $\omega$, when the height
$H$ of the excursion is high. Let us give a precise statement. Define
the critical height
\[
\label{hTdef} \hP_t\defeq\log t - \log\log t,\qquad t \ge
\ee^{\ee}.\vadjust{\goodbreak}
\]

Heuristics suggest (and it would follow from later results) that when
$H>\hP_t$, $\tau$ is on the order of $\ee^H>\frac{t}{\log t}$.
Proposition~\ref{propexitonevalley} shows that the distance between
$\tau$ and $E_\omega[\tau]\eg$ (for a suitable coupling) is no larger
than $t^\beta\ll\frac{t}{\log t}$ in quenched average when $H>\hP_t$,
in agreement with our aim.

%
\begin{proposition}\label{propexitonevalley}
We may enlarge the probability space in order to introduce an
exponential random variable $\eg$ of parameter $1$, independent of
$\omega$, such that, for some $\beta<1$, as $t\to\infty$,
%
%
\begin{equation}
\label{eqnonevalley} \Pp\bigl(E_\omega\bigl[ \bigl|\tau-E_\omega[
\tau]\eg\bigr| \bigr]>t^\beta,H\geq\hP_t \bigr)=o
\bigl(t^{-\kappa}\bigr).
\end{equation}
\end{proposition}

This proposition can equivalently be phrased, using \eqref
{iglehartthm}, as
\[
\Pp\bigl(W^1_\omega\bigl(\tau,E_\omega[\tau]\eg
\bigr)>t^\beta|H\geq\hP_t \bigr)=o \biggl(\frac1{\log t}
\biggr),
\]
where $\eg$ is an exponential random variable of parameter $1$
independent of~$\omega$.

\subsection{``Good'' environments}
The proof relies on a precise control of the geometry of a typical
valley, namely, that it is not too wide and smooth enough.
Let us define the maximal ``increments'' of the potential in a window
$ [x,y ]$ by
\begin{eqnarray*}
\label{incremupgene} V^{\uparrow}(x,y)&:=& \max_{x \le u \le v \le y}
\bigl(V(v)-V(u)\bigr),\qquad x<y,
\\
\label{incremdowngene} V^{\downarrow}(x,y)&:=& \min_{x \le u \le v
\le y}
\bigl(V(v)-V(u)\bigr),\qquad x<y.
\end{eqnarray*}

Then, we introduce the following events:
\begin{eqnarray*}
\Omega^{(1)}_t &\defeq&\{ e_1 \le C \log t
\},
\\
\Omega^{(2)}_t &\defeq&\bigl\{ \max\bigl\{
-V^{\downarrow}(0,T_H), V^{\uparrow}(T_H,e_1)
\bigr\} \le\alpha\log t \bigr\},
\\
\Omega^{(3)}_t &\defeq&\bigl\{ R^- \le(\log
t)^4 t^\alpha\bigr\},
\end{eqnarray*}
where $\max\{0,1-\kappa\}<\alpha<\min\{1,2-\kappa\}$ is arbitrary, and
$R^-$ is defined by
\[
R^- := \sum_{x=-\infty}^{-1} \Biggl(1+2 \sum
_{y=x+2}^{0} \ee^{V(y)-V(x+1)} \Biggr)
\Biggl( \ee^{-V(x+1)} +2 \sum_{y=-\infty
}^{x-1}
\ee^{-V(y+1)} \Biggr).
\]

We define the set of ``good'' environments at time $t$ by
%
%
\begin{equation}
\label{eqgoodenv} \Omega_t\defeq\Omega^{(1)}_t
\cap\Omega^{(2)}_t \cap\Omega^{(3)}_t.
\end{equation}
By the following result, ``good'' environments are asymptotically
typical on $\{H\geq\hP_t \}$.
%
%
\begin{lemma}
\label{lemgoodenv}
The event $\Omega_t$ satisfies
\[
\label{eqestimatevalley2} \Pp\bigl( \Omega_t^c, H \ge
\hP_t\bigr)=o\bigl(t^{-\kappa}\bigr),\qquad t \to\infty.\vadjust{\goodbreak}
\]
\end{lemma}

The proof of this result is easy but technical and therefore postponed
to Appendix~\ref{subsecprooflemma-singlevalley}.

\subsection{Preliminary results}
\label{subsechproc}

In order to finely estimate the time spent in a deep valley, we
decompose the passage from $0$ to $e_1$ into the sum of a random
geometrically distributed number, denoted by $N$, of unsuccessful
attempts to reach $e_1$ from $0$ (i.e., excursions of the particle from
$0$ to $0$ which do not hit $e_1$), followed by a successful attempt.
More precisely, $N$ is a geometrically distributed random variable with
parameter $1-p$
satisfying
%
%
\begin{equation}
\label{1-p} 1-p = \frac{\omega_0}{\sum_{x=0}^{e_1-1}
\ee^{V(x)}}= \frac{\omega_0}{M_2 \ee^H},
\end{equation}
and we can write $ \tau(e_1)=\sum_{i=1}^{N} F_i +G,$ where the $F_i$'s
are the durations of the successive i.i.d. failures and $G$ that of the
first success. The accurate estimation of the time spent by each
(successful and unsuccessful) attempt leads us to consider two
$h$-processes where the random walker evolves in two modified
potentials, one corresponding to the conditioning on a failure (this
potential is denoted by $\widehat{V}$ in~\cite{limitlaws}, page 2494)
and the other to the conditioning on a success (denoted by $\bar{V}$ in
\cite{limitlaws}, page 2497). Note that this approach was first
introduced by three of the authors in~\cite{limitlaws} to estimate the
quenched Laplace transform of the occupation time of a deep valley in
the case $0<\kappa<1$. We refer to~\cite{limitlaws} for more details on
these two $h$-processes. Moreover, using the properties of ``good''
environments introduced above, we can prove the following useful
lemmata, whose proofs are postponed to Appendix~\ref
{subsecprooflemma-singlevalley}.
%
%
\begin{lemma}
\label{lemmafailurebound}
For all $t\ge1,$ we have on $\Omega_t$,
%
%
\begin{eqnarray}
\label{eqboundVarF} \Var_\omega(F) &\le&C (\log t)^4
t^\alpha,
\\
M_2 &\le& C \log t,\label{eqboundM2}
\\
\vert\widehat M_1 -M_1 \vert&\le& o\bigl(t^{-\delta}
\bigr)M_1,\label{eqboundM1}
\end{eqnarray}
with $\delta\in(0,1-\alpha)$ and where $\widehat M_1$ is defined by
the relation $ E_{\omega} [ F ] = 2 \omega_0
\widehat M_1.$
\end{lemma}

%
\begin{lemma}
\label{lemmafsuccessbound}
For all $t\ge1,$ we have on $\Omega_t$,
\[
E_{\omega}[G] \le C (\log t)^4 t^\alpha.
\]
\end{lemma}

\subsection{Definition of the coupling} \label{subseccoupling}
We recall here the coupling from~\cite{aging} between the quenched
distribution of the random walk before time $\tau$ and an exponential
random variable $\eg$ of parameter $1$ independent of $\omega$. Given
$\omega$ and $\eg$, let us define
\[
N\defeq\biggl\lfloor-\frac{1}{\log(1-p(\omega))}\eg\biggr
\rfloor,
\]
where $p(\omega)=P_{0,\omega}(\tau(0)<\tau(e_1))$; cf. \eqref{1-p}.
Note that, conditionally on $\omega$, $N$ is a geometric random
variable of parameter $1-p$, just like the number of returns to~$0$
before the walk reaches $e_1$.

Given $\omega$ and $\eg$ (and hence, $N$), the random walk is sampled
as usual as a Markov chain, except that the number of returns to $0$ is
conditioned on being equal to $N$, which amounts to saying that when
the walk reaches $0$ for the first $N$ times, it is conditioned on
coming back to $0$ before reaching $e_1$ (this is still a Markov chain,
namely, the $h$-process associated to $\widehat{V}$; see
Section~\ref{subsechproc}), while on the $(N+1)$th visit of $0$ it is
conditioned
on reaching $e_1$ first (this is the $h$-process associated to~$\bar V$).
Due to the definition of $N$, the distribution of the walk given
$\omega$ only is $P_{0,\omega}$.

\subsection{\texorpdfstring{Proof of Proposition \protect\ref{propexitonevalley}}{Proof of Proposition 1}}\label{subseconevalley}

We consider the same decomposition as in Section~\ref{subsechproc},
that is,
$\tau=F_1+\cdots+F_N+G.$ By Wald identity, $E_\omega[\tau]=E_\omega
[N]E_\omega[F]+E_\omega[G]$. Thus, we have
\begin{eqnarray*}
\bigl|\tau-E_\omega[\tau]\eg\bigr| & \leq&\bigl|F_1+\cdots+F_N-N
E_\omega[F] \bigr|+E_\omega[F] \bigl|N-E_\omega[N]\eg\bigr|
\\
&&{} + G+E_\omega[G]\eg.
\end{eqnarray*}
Let us consider each term, starting with the last two (with same
$P_\omega$-expecta\-tion). If we choose $\beta$ such that $\alpha
<\beta
<1$, then by Lemmas~\ref{lemmafsuccessbound} and~\ref{lemgoodenv} we
have, for large $t$,
%
%
\begin{equation}
\Pp\biggl(E_\omega[G]\geq\frac{t^\beta}4, H\geq\hP_t
\biggr) \leq\Pp\bigl((\Omega_t)^c, H\geq
\hP_t \bigr) = o \bigl(t^{-\kappa
} \bigr).\label{emajorg}
\end{equation}
We turn to the first one. Conditioning first on $N$ [which is
independent of $(F_i)_i$] and then applying the Cauchy--Schwarz inequality,
we have
\begin{eqnarray*}
E_\omega\bigl[ \bigl|F_1+\cdots+F_N-N
E_\omega[F] \bigr| \bigr] &\leq& E_\omega\bigl[\Var_\omega
(F_1+\cdots+F_N |N )^{1/2} \bigr]
\\
& = &E_\omega\bigl[N^{1/2} \bigr]\Var_\omega(F)^{1/2}.
\end{eqnarray*}
Furthermore, $E_\omega[N^{1/2}]\leq E_\omega
[N]^{1/2}=((1-p)^{-1}-1)^{1/2}\leq(M_2)^{1/2}\ee^{H/2}\omega_0^{-1/2}$
and $\omega_0\geq\frac12,$ $\Pp$-almost surely. Thus, using
Lemma~\ref{lemmafailurebound} to bound $\Var_\omega(F)$, we get
\begin{eqnarray*}
&&\Pp\biggl(E_\omega\bigl[\bigl|F_1+\cdots+F_N-N
E_\omega[F] \bigr| \bigr]>\frac{t^\beta}4,H\geq\hP_t
\biggr)
\\
&&\qquad\leq\Pp\bigl((\Omega_t)^c,H\geq\hP_t
\bigr) + \Pp\biggl((M_2)^{1/2}\ee^{ H/2}\geq
\frac{t^{\beta-
\alpha/2}}{C(\log t)^2},H\geq\hP_t \biggr).
\end{eqnarray*}
As before, the first term is $o(t^{-\kappa})$. And the second one is
less than
\begin{eqnarray*}
&&P\bigl(M_2\geq(\log t)^2, H\geq\hP_t
\bigr)+P \biggl(\ee^{ H/2}\geq\frac
{t^{\beta-\alpha/2}}{C(\log t)^3} \biggr)
\\
&&\qquad\leq\frac{P(H\geq\hP_t)}{(\log t)^2}E[M_2|H\geq\hP_t]+P
\biggl(
\ee^H\geq\frac{t^{2\beta-\alpha}}{C^2(\log t)^6} \biggr).
\end{eqnarray*}
Each term is $o(t^{-\kappa})$ if we additionally impose $\frac
{1+\alpha
}2<\beta<1$, due to \eqref{eqnm1m2sup} and~\eqref{iglehartthm}.

Finally, we have
\begin{eqnarray*}
\bigl|N-E_\omega[N]\eg\bigr| & =& \biggl| \biggl\lfloor\frac
{1}{-\log(1-p)}\eg
\biggr
\rfloor- \biggl(\frac
{1}{p}-1 \biggr)\eg\biggr|
\\
& \leq&\biggl(1+ \biggl|-\frac{1}{\log(1-p)}-\frac{1}{p} \biggr|
\biggr)\eg+1,
\end{eqnarray*}
and the function $x\mapsto-\tfrac{1}{\log(1-x)}-\tfrac{1}{x}$ extends
continuously on $[0,1]$ and is thus bounded by a constant $C$, hence,
\begin{eqnarray*}
&&\Pp\biggl(E_\omega[F]E_\omega\bigl[ \bigl|N-E_\omega[N]
\eg\bigr| \bigr]\geq\frac{t^\beta}4, H\geq\bP_t \biggr)
\\
&&\qquad\leq\Pp\biggl(E_\omega[F]\geq\frac{t^\beta}{4C}, \
H\geq
\hP_t \biggr)
\\
&&\qquad\leq\Pp\bigl((\Omega_t)^c,H\geq\hP_t
\bigr)+\frac{8C P(H\geq\hP_t)}{t^\beta}\Ep[M_1|H\geq\hP_t]
\end{eqnarray*}
for large $t$, due to \eqref{eqboundM1}, recalling that $ E_{\omega}
[ F ] = 2 \omega_0 \widehat M_1$ (see Lemma \ref
{lemmafailurebound}).
We conclude as before that this is negligible compared to $t^{-\kappa}$.

Therefore, gathering all these estimates gives Proposition \ref
{propexitonevalley}.

\section{Independence of the deep valleys} \label{seciidvalleys}

The independence between deep valleys goes through imposing these
valleys to be disjoint (i.e., $a_i>d_{i-1}$ for all $i$) and neglecting
the time spent on the left of a valley while it is being crossed
(i.e., the time spent on the left of $a_i$ before $d_i$ is reached).

\textit{NB}. All the results and proofs from this section hold for
any parameter $\kappa>0$.

For any integers $x,y,z$, let us define
\[
\tautx{z}(x,y)\defeq\#\bigl\{\tau(x)\leq k\leq\tau(y)\dvtx
X_k\leq
z\bigr\},
\]
the time spent on the left of $z$ between the first visit to $x$ and
the first visit to $y$, and
\[
\tautx{z}\defeq\#\bigl\{k\geq\tau(z)\dvtx X_k\leq z\bigr\},
\]
the total time spent on the left of $z$ after the first visit to $z$.

Let us consider the event
\[
\label{defnonoverlap} \nonoverlap\defeq\{0<a_1\}\cap\bigcap
_{i=1}^{K_n-1}\{d_i<a_{i+1}\},
\]
which means that the large valleys before $e_n$ lie entirely on $\Z_+$
and do not overlap.
The following two propositions will enable us to reduce to i.i.d. deep valleys.

%
\begin{proposition}\label{propoverlap}
We have
\[
\label{eqlimtauoverlap} P\bigl(\nonoverlap\bigr)\limites{} {n} 1.
\]
\end{proposition}

\begin{pf}
Choose $\eps>0$ and define the event
\[
A_K(n)\defeq\bigl\{K_n\leq(1+\eps)C_I(
\log n)^\kappa\bigr\}.
\]
Since $K_n$ is a binomial random variable of mean $nq_n\sim_n C_I(\log
n)^\kappa$, it follows from the law of large numbers that $P(A_K(n))$
converges to $1$ as $n\to\infty$. On the other hand, if the event
$\overlap$ occurs, then there exists $1\leq i\leq K_n$ such that there
is at least one high excursion among the first $D_n$ excursions to the
right of $d_{i-1}$ (with $d_0=0$). Thus,
\begin{eqnarray*}
P\bigl(\overlap\bigr) & \leq& P\bigl(A_K(n)^c\bigr)+(1+
\eps)C_I(\log n)^\kappa\bigl(1-(1-q_n)^{D_n}
\bigr)
\\
& \leq& o(1)+(1+\eps)C_I(\log n)^\kappa q_n
D_n = o(1).
\end{eqnarray*}
Indeed, for any $0<u<1$ and $\alpha>0$, we have $1-(1-u)^\alpha\leq
\alpha u$ by concavity of $u\mapsto1-(1-u)^\alpha$.
\end{pf}

%
\begin{proposition}\label{proplimtaut}
Under $\Pp$,
\[
\frac{1}{n^{1/\kappa}}\sum_{i=1}^{K_n}E_\omega
\bigl[\tautx{a_i}(b_i,d_i)\bigr]=
\frac{1}{n^{1/\kappa}}\sum_{p=0}^{n-1}E_\omega
\bigl[\tautx{e_{p-D_n}}(e_p,e_{p+1})
\bigr]\indic{H_p\geq h_n}\limites{(p)} {n}0.
\]
\end{proposition}

\begin{pf}
The equality is trivial from the definitions. Note that the terms in
the second expression have the same distribution under $\Pp$ because of
Lemma~\ref{lemenegatifs}. As $E_\omega[\tautx
{e_{-D_n}}(0,e_1)]\indic
{H\geq h_n}$ is not integrable for $0<\kappa\leq1$, we introduce the event
\begin{eqnarray*}
A_n &\defeq&\bigl\{\mbox{for }i=1,\ldots,K_n,
H_{\sigma(i)}\leq V(a_i)-V(b_i)\bigr\}
\\
&\hphantom{:}=& \bigcap_{p=0}^{n-1}
\{H_p<h_n\}\cup\bigl\{h_n\leq
H_p\leq V(e_{p-D_n})-V(e_p)\bigr\}.
\end{eqnarray*}
Let us prove that our choice of $D_n$ ensures $\Pp((A_n)^c)=o_n(1)$. By
Lemma~\ref{lemenegatifs}, we have
$
\Pp((A_n)^c)
\leq n \Pp(H\geq h_n, H> V(e_{-D_n})).
$
Then, let us choose $0<\gamma'<\gamma''<\gamma$ [cf. \eqref
{eqndefdn}] and define $l_n\defeq\frac{1+\gamma'}{\kappa}\log
n$. We get
%
%
\begin{equation}
\Pp\bigl((A_n)^c\bigr) \leq n \bigl( P(H\geq
l_n)+P(H\geq h_n)\Pp\bigl(V(e_{-D_n})<l_n
\bigr) \bigr). \label{eqn4697}
\end{equation}
Equation \eqref{iglehartthm} gives $P(H\geq l_n)\sim_n C_I\ee
^{-\kappa
l_n}=C_I n^{-(1+\gamma')}$, and $P(H\geq h_n)$ $\sim_n C_I
n^{-1}(\log n)^\kappa$. Under $\Pp$, $V(e_{-D_n})$ is the sum of $D_n$
i.i.d. random\vadjust{\goodbreak} variables distributed like $-V(e_1)$. Therefore, for any
$\lambda>0$,
\[
\Pp\bigl(V(e_{-D_n})<l_n\bigr) \leq\ee^{\lambda l_n} E
\bigl[\ee^{-\lambda(-V(e_1))}\bigr]^{D_n}.
\]
Since $\frac{1}{\lambda}\log E[\ee^{-\lambda(-V(e_1))}]\to
-E[-V(e_1)]\in[-\infty,0)$ as $\lambda\to0^+$, we can choose
$\lambda
>0$ such that $\log E[\ee^{-\lambda(-V(e_1))}]<-\lambda A\frac
{1+\gamma''}{1+\gamma}$ [where $A$ was defined after \eqref{eqndefdn}],
hence, $E[\ee^{-\lambda(-V(e_1))}]^{D_n}\leq n^{-\lambda
{(1+\gamma
'')}/{\kappa}}$. This gives the bound\break $\Pp(V(e_{-D_n})<l_n)\leq
n^{-\lambda{(\gamma''-\gamma')}/{\kappa}}$. Using these
estimates in
\eqref{eqn4697} concludes the proof that $\Pp((A_n)^c)=o_n(1)$.

Let us now prove the proposition itself. By the Markov inequality, for all
$\delta>0$,
%
%
\begin{eqnarray}\label{eqn1564}
&& \Pp\Biggl(\frac{1}{n^{1/\kappa}}\sum_{p=0}^{n-1}E_\omega
\bigl[\tautx{e_{p-D_n}}(e_p,e_{p+1})\bigr]
\indic{H_p\geq h_n}>\delta\Biggr)
\nonumber
\\
&&\qquad\leq\Pp\bigl((A_n)^c\bigr)+\frac{1}{\delta n^{1/\kappa
}}\Ep
\Biggl[\sum_{p=0}^{n-1}E_\omega
\bigl[\tautx{e_{p-D_n}}(e_p,e_{p+1})\bigr]
\indic{H_p\geq h_n}\ind_{A_n} \Biggr]
\\
&&\qquad\leq o_n(1) + \frac{n}{\delta n^{1/\kappa}}\Ep\bigl
[E_\omega
\bigl[\tautx{e_{-D_n}}(0,e_1)\bigr]\indic{H\ge
h_n, H<V(e_{-D_n})} \bigr].\nonumber
\end{eqnarray}
Note that we have $E_\omega[\taut^{(e_{-D_n})}(0,e_1)]=E_\omega
[N]E_\omega[T_1]$, where $N$ is the number of crossings from
$e_{-D_n}+1$ to $e_{-D_n}$ before the first visit at $e_1$, and $T_1$
is the time for the random walk to go from $e_{-D_n}$ to $e_{-D_n}+1$
(for the first time, e.g.); furthermore, these two terms are
independent under $\Pp$. Using \eqref{eqnzeitounip}, we have
\begin{eqnarray*}
E_\omega[N]&=&\frac{P_{0,\omega}(\tau(e_{-D_n})<\tau
(e_1))}{P_{e_{-D_n}+1,\omega}(\tau(e_1)<\tau(e_{-D_n}))}= \sum
_{0\leq x<e_1}
\ee^{V(x)-V(e_{-D_n})}\\
&=& M_2\ee^{H-V(e_{-D_n})},
\end{eqnarray*}
hence, on the event $\{H<V(e_{-D_n})\}$, $E_\omega[N]\leq M_2$.

The length of an excursion to the left of $e_{-D_n}$ is computed as
follows, due to~\eqref{eqnzeitounie3}:
\[
E_\omega[T_1]=E_{e_{-D_n},\omega}\bigl[\tau(e_{-D_n}+1)
\bigr]\leq2\sum_{x\leq
e_{-D_n}}\ee^{-(V(x)-V(e_{-D_n}))}.
\]
The law of $(V(x)-V(e_{-D_n}))_{x\leq e_{-D_n}}$ under $\Pp$ is $\Pp$
because of Lemma~\ref{lemenegatifs}. Therefore,
\[
\Ep\bigl[E_\omega[T_1]\bigr]\leq2\Ep\biggl[\sum
_{x\leq0}\ee^{-V(x)} \biggr]=2\Ep[R_-]<\infty,
\]
with \eqref{eqnr-}. Then, we conclude that the right-hand side of
\eqref{eqn1564} is less than
$
o_n(1)+2 \delta^{-1 }n^{1-1/\kappa}\Ep[R_-]E[M_2\indic{H\ge h_n}].
$
Since Lemma~\ref{lemrenewalestimates}\vspace*{1pt} gives the bound\vadjust{\goodbreak}
$E[M_2\indic
{H\ge h_n}]\leq CP(H\ge h_n)\sim_n C'\ee^{-\kappa h_n}=C'n^{-1}(\log
n)^\kappa$, this whole expression converges to 0, which concludes the
proof of the proposition.
\end{pf}

\section{Fluctuation of interarrival times} \label{secinterarrival}

For any $x\leq y$, recall that the interarrival time $\tau(x,y)$
between sites $x$ and $y$ is defined by
$
\tau(x,y) \defeq\inf\{ n \ge0\dvtx X_{\tau(x)+n}=y \}.
$
Then, let
\[
\tau_{\mathrm{IA}}\defeq\sum_{i=0}^{K_n}
\tau(d_i,b_{i+1}\wedge e_n)=\sum
_{p=0}^{n-1}\tau(e_p,e_{p+1})
\indic{H_p<h_n}
\]
(with $d_0=0$) be the time spent at crossing small excursions before
$\tau(e_n)$. The aim of this section is the following bound on the
fluctuations of $\tau_{\mathrm{IA}}$.
%
%
\begin{proposition}\label{proptauia}
For any $0<\kappa<2$, under $\Pp$,
%
%
\begin{equation}
\frac{1}{n^{1/\kappa}}E_\omega\bigl[ \bigl|\tau_{\mathrm
{IA}}-E_\omega[\tau_{\mathrm{IA}}] \bigr|
\bigr]\limites{(p)} {n}0.\label{eqnlimitefluctutia}
\end{equation}
\end{proposition}

This proposition holds for $0<\kappa<1$ in a simple way: we have, in
this case, using Lemmas~\ref{lemenegatifs} and~\ref{lemboundEe},
\begin{eqnarray*}
\Ep\bigl[E_\omega[\tau_{\mathrm{IA}}]\bigr]&=&n\Ep\bigl
[E_\omega\bigl[
\tau(e_1)\bigr]\indic{H< h_n}\bigr] \leq n\Ep
\bigl[2M'_1M_2\ee^H\indic{H<
h_n}\bigr]
\\
&\leq& C n\ee^{(1-\kappa)h_n}=o\bigl(n^{1/\kappa}\bigr),
\end{eqnarray*}
hence, $n^{-1/\kappa}E_\omega[\tau_{\mathrm{IA}}]$ itself converges
to 0 in
$L^1(\Pp
)$-norm and thus in probability.

We now consider the case $1\leq\kappa<2$. The proposition will directly
follow from the fact that, under $\Pp$,
\[
\frac{1}{n^{2/\kappa}}\Var_\omega(\tau_{\mathrm{IA}})\limites
{(p)} {n}0,
\]
which in turn will come from Lemma~\ref{lemvartia} proving $\Ep
[\Var_\omega(\tau_{\mathrm{IA}})]=o(n^{2/\kappa})$. However, a
specific caution is
necessary in the case $\kappa=1$; indeed, $\Var_\omega(\tau
_{\mathrm{IA}})$ is not
integrable in this case, because of the rare but significant
fluctuations originating from the time spent by the walk when it
backtracks into deep valleys. Our proof in this case consists of first
proving that we may neglect in probability (using a first-moment
method) the time spent backtracking into these deep valleys; and then
that this brings us to the computation of the variance of $\tau
_{\mathrm{IA}}$ in an
environment where small excursions have been substituted for the high
ones (thus removing the nonintegrability problem).\looseness=1

Section~\ref{subsecreducsmall} is dedicated to this reduction to
an integrable setting, which is only involved in the case $\kappa=1$ of
Proposition~\ref{proptauia} and of the theorems (but holds in greater
generality), while Section~\ref{subsecvartia} states and proves
the bounds on the variance, implying Proposition~\ref{proptauia}.

\subsection{\texorpdfstring{Reduction to small excursions $($required for the case $\kappa=1)$}
{Reduction to small excursions (required for the case kappa=1)}}\label{subsecreducsmall}

Let \mbox{$h>0$}. Let us denote by $d_-$ the right end of the first excursion
on the left of 0 that is higher than $h$:
\[
d_- \defeq\max\{e_p\dvtx p\leq0, H_{p-1}\geq h\}.
\]
Remember $\tautx{d_-}(0,e_1)$ is the time spent on the left of $d_-$
before the walk reaches~$e_1$.

%
\begin{lemma}\label{lemremonteedifficile}
There exists $C>0$, independent of $h$, such that
%
%
\begin{equation}
\label{eqntimeleftd} \bEp\bigl[\tautx{d_-}(0,e_1)\indic{H< h}\bigr]
\leq C \cases{ %
\ee^{-(2\kappa-1)h},&\quad$\mbox{if }\kappa<1,$
\vspace*{2pt}\cr
h\ee^{- h},&\quad$\mbox{if }\kappa=1,$
\vspace*{2pt}\cr
\ee^{-\kappa h},&\quad$\mbox{if }\kappa>1.$}
\end{equation}
\end{lemma}

\begin{pf}
Let us decompose $\tautx{d_-}(0,e_1)$ into the successive excursions to
the left of $d_-$:
\[
\tautx{d_-}(0,e_1)=\sum_{m=1}^N
T_m,
\]
where $N$ is the number of crossings from $d_-+1$ to $d_-$ before $\tau
(e_1)$, and $T_m$ is the time for the walk to go from $d_-$ to $d_-+1$
on the $m$th time. Under $P_\omega$, the times $T_m$, $m\geq1$, are
i.i.d. and independent of $N$ [i.e., more properly, the sequence
$(T_m)_{1\leq m\leq N}$ can be prolonged to an infinite sequence with
these properties]. We have, using the Markov property and then \eqref
{eqnzeitounip},
\[
E_\omega[N]=\frac{P_{0,\omega}(\tau(d_-)<\tau
(e_1))}{P_{d_-+1,\omega
}(\tau(e_1)<\tau(d_-))}=\sum_{0\leq x<e_1}
\ee^{V(x)-V(d_-)}
\]
and, from \eqref{eqnzeitounie3},
$
E_\omega[T_1]=E_{d_-,\omega}[\tau(d_-+1)]\leq2\sum_{x\leq d_-}\ee
^{-(V(x)-V(d_-))}.
$
Therefore, by Wald identity and Lemma~\ref{lemenegatifs},
%
%
\begin{eqnarray}\label{eqn7878}
&&\bEp\bigl[\tautx{d_-}(0,e_1)\indic{H< h}\bigr]\nonumber\\
&&\qquad= \Ep
\bigl[E_\omega[N]E_\omega[T_1]\indic{H< h}\bigr]
\\
&&\qquad\leq2 E \biggl[\sum_{0\leq x<e_1}\ee^{V(x)}
\indic{H< h} \biggr] \Ep\bigl[\ee^{-V(d_-)}\bigr] E \biggl[\sum
_{x\leq0}\ee^{-V(x)} \Big|\Lambda(h) \biggr],\nonumber
\end{eqnarray}
where $\Lambda(h)\defeq\{\forall x\leq0, V(x)\geq0\}\cap\{
H_{-1}\geq
h\}$. The first expectation can be written as $E[M_2\ee^H \indic{H<
h}]$. For the second one, note that $d_-=e_{-W}$, where $W$ is a
geometric random variable of parameter $q\defeq P(H\geq h)$; and,
conditional on $\{W=n\}$, the distribution of $(V(x))_{e_{-W}\leq x\leq
0}$ under $\Pp$ is the same as\vadjust{\goodbreak} that of $(V(x))_{e_{-n}\leq x\leq0}$
under $\Pp(\cdot|\mbox{for }p=0,\ldots,n-1, H_{-p}< h)$. Therefore,
\begin{eqnarray*}
\Ep\bigl[\ee^{-V(d_-)}\bigr] = \Ep\bigl[E\bigl[\ee^{-V(e_1)}|H<h
\bigr]^W \bigr] = \frac{q}{1-(1-q)E[\ee^{V(e_1)}|H<h]},
\end{eqnarray*}
and $(1-q)E[\ee^{V(e_1)}|H<h]$ converges to $E[\ee^{V(e_1)}]<1$ when
$h\to\infty$ [the inequality comes from assumption (b)], hence,
this quantity is uniformly bounded from above by $c<1$ for large $h$.
In addition, \eqref{iglehartthm} gives $q\sim C_I\ee^{-\kappa h}$ when
$h\to\infty$, hence,
$
\Ep[\ee^{-V(d_-)}]\leq C\ee^{-\kappa h},
$
where $C$ is independent of $h$.

Finally, let us consider the last term of \eqref{eqn7878}. We have
\begin{eqnarray*}
&& E \biggl[\sum_{x\leq0}\ee^{-V(x)} \Big|\Lambda
(h) \biggr]
\\
&&\qquad= E \biggl[\sum_{e_{-1}<x\leq0}\ee^{-V(x)}
\Big|H_{-1}\geq h \biggr] +\Ep\biggl[\sum_{x\leq e_{-1}}
\ee^{-(V(x)-V(e_{-1}))} \biggr]E\bigl[\ee^{-V(e_{-1})}\bigr]
\\
& &\qquad\leq E\bigl[M'_1|H\geq h\bigr]+\Ep[R_-]E\bigl[
\ee^{V(e_1)}\bigr],
\end{eqnarray*}
hence, using Lemma~\ref{lemrenewalestimates}, \eqref{eqnr-} and
$V(e_1)\leq0$, this term is bounded by a constant. The statement of
the lemma then follows from the application of Lemma \ref
{lemrenewalestimates} to the expectation $E[M_2\ee^H\indic{H< h}]$.
\end{pf}

The part of the interarrival time $\tau_{\mathrm{IA}}$ spent at
backtracking in high
excursions can be written as follows:
\begin{eqnarray*}
\widetilde{\tau}_{\mathrm{IA}} &\defeq&\tautx{d_-}(0,b_1\wedge e_n)+
\sum_{i=1}^{K_n}\tautx{d_i}(d_i,b_{i+1}
\wedge e_n)
\\
&\hphantom{}=&\sum_{p=0}^{n-1} \tautx{d(e_p)}(e_p,e_{p+1})
\indic{H_p<h_n},
\end{eqnarray*}
where, for $x\in\Z$,
$
d(x)\defeq\max\{e_p\dvtx p\in\Z, e_p\leq x, H_{p-1}\geq h_n\}.
$
In particular, $d(0)=d_-$ in the previous notation with $h=h_n$.

Note that under $\bPp$, because of Lemma~\ref{lemenegatifs}, the
terms of the above sum have the same distribution as $\tautx
{d(0)}(0,e_1)\indic{H< h_n}$, hence,
\[
\bEp[\widetilde{\tau}_{\mathrm{IA}}]= n\bEp\bigl[\tautx
{d(0)}(0,e_1)\indic{H<
h_n}\bigr].
\]
Thus, for $\bEp[\widetilde{\tau}_{\mathrm{IA}}]$ to be negligible
with respect to
$n^{1/\kappa}$, it suffices that the expectation on the right-hand side
be negligible with respect to $n^{1/\kappa-1}$. In particular, for
$\kappa=1$, it suffices that it converges to 0, which is readily seen
from \eqref{eqntimeleftd}. Thus, for $\kappa=1$,
%
%
\begin{equation}
\label{eqnlimtiatilde} n^{-1/\kappa} \bEp[\widetilde{\tau}_{\mathrm{IA}
}]\limites{} {n}0,
\end{equation}
hence, in particular, $n^{-1/\kappa}E_\omega[\widetilde{\tau}_{\mathrm{IA}}]\to
0$ in
probability under $\bPp$. Note that \eqref{eqnlimtiatilde} actually
holds for any $\kappa\geq1$.\vadjust{\goodbreak}

Let us introduce the modified environment, where independent small
excursions are substituted for the high excursions. In order to avoid
obfuscating the redaction, we will only introduce little notation
regarding this new environment.

Let us enlarge the probability space in order to accommodate a new
family of independent excursions indexed by $\N^*\times\Z$ such that
for all $n,k$ the excursion with index $(n,k)$ has the same
distribution as $(V(x))_{0\leq x\leq e_1}$ under $P( \cdot|H<
h_n)$. Thus we are given, for every $n\in\N^*$, a countable family of
independent excursions lower than $h_n$. For every fixed $n$, we define
the \textit{modified environment of height less than $h_n$} by replacing
all the excursions of $V$ that are higher than $h_n$ by new independent
ones that are lower than $h_n$. Because of Lemma~\ref{lemenegatifs},
this construction is especially natural under $\Pp$, where it has
stationarity properties.

In the following, we will denote by $P'$ the law of the modified
environment relative to the height $h_n$ given in the context [hence,
also a definition of $(\Pp)'$, e.g.].

\begin{remark*}\label{remer}
Repeating the proof done under $\Pp$ for $(\Pp)'$, we see that $R_-$
still has all finite moments in the modified environment, and
that\label{rem:e'r}
these moments are bounded uniformly in $n$. In particular, the bound
for the quantity $\Ep[(M'_1)^\alpha(M_2)^\beta\ee^{\gamma H}\indic{H<
h_n}]$ given in Lemma~\ref{lemrenewalestimates} is unchanged for
$(\Ep)'$ $[$writing $M'_1=R_-+\sum_{0\leq x<e_1}\ee^{-V(x)}$ and using
$(a+b)^\alpha\leq2^\alpha(a^\alpha+b^\alpha)]$. On the other hand,
\begin{eqnarray*}
E'[R] &=& \sum_{p=0}^\infty
E'\bigl[\ee^{V(e_p)}\bigr]E' \biggl[\sum
_{e_p\leq
x<e_{p+1}}\ee^{V(x)-V(e_p)} \biggr]
\\
&= &\sum_{p=0}^\infty E\bigl[
\ee^{V(e_1)}|H< h_n\bigr]^p E
\bigl[M_2\ee^H|H< h_n\bigr],
\end{eqnarray*}
and $E[\ee^{V(e_1)}|H< h_n]\leq c$ for some $c<1$ independent of $n$
because this expectation is smaller than $1$ for all $n$ and it
converges toward $E[\ee^{V(e_1)}]<1$ as $n\to\infty$. Hence, by Lemma
\ref{lemrenewalestimates},
%
%
\begin{equation}
\label{eqner} \mbox{if $\kappa=1$}\qquad E'[R] \leq C h_n.
\end{equation}
This is the only difference that will appear in the following computations.
\end{remark*}

Assuming that $d(0)$ keeps being defined with respect to the usual
heights, \eqref{eqntimeleftd} (with $h=h_n$) is still true for the
walk in the modified environment. Indeed, the change only affects the
environment on the left of $d(0)$, hence, the only difference in the
proof involves the times $T_m$; in \eqref{eqn7878}, one should
substitute $(\Ep)'$ for $E[\cdot|{\Lambda(h)}]$, and this factor is
uniformly bounded in both cases because of the above remark about $R_-$.

We deduce that the time $\widetilde{\tau}_{\mathrm{IA}}'$, defined
as similar to
$\widetilde{\tau}_{\mathrm{IA}}$ except that the excursions on the
left of the points
$d(e_i)$ (i.e., the times similar to $T_m$ in the previous proof) are
performed in the modified environment, still satisfies, for $\kappa=1$,
%
%
\begin{equation}
\label{eqnlimtiatildee}
n^{-1/\kappa} \bEp\bigl[\widetilde{
\tau}_{\mathrm{IA}}'\bigr]\limites{} {n}0.
\end{equation}

Now note that
%
%
\begin{equation}
\label{eqtiaprime} \tau_{\mathrm{IA}}'\defeq\tau_{\mathrm
{IA}}-\widetilde{\tau}_{\mathrm{IA}}+\widetilde{
\tau}_{\mathrm{IA}}'
\end{equation}
is the time spent at crossing the (original) small excursions, in the
environment where the high excursions have been replaced by new
independent small excursions. Indeed, the high excursions are only
involved in $\tau_{\mathrm{IA}}$ during the backtracking of the walk
to the left of
$d(e_i)$ for some $0\leq i< n$. Assembling \eqref{eqnlimtiatilde}
and \eqref{eqnlimtiatilde}, it is equivalent (for $\kappa=1$) to
prove \eqref{eqnlimitefluctutia} or
\[
n^{-1/\kappa} E_\omega\bigl[ \bigl|\tau_{\mathrm{IA}}'-E_\omega
\bigl[\tau_{\mathrm{IA}}'\bigr] \bigr| \bigr]\limites{(p)} {n}0,
\]
and it is thus sufficient to prove
$
\Ep[\Var_\omega(\tau_{\mathrm{IA}}')]=o_n(n^{2/\kappa}).
$

\subsection{\texorpdfstring{Bounding the variance of $\tau_{\mathrm{IA}}$}
{Bounding the variance of tauIA}} \label{subsecvartia}
Because of the previous subsection, Proposition~\ref{proptauia} will
follow from the next lemma.

%
\begin{lemma}\label{lemvartia}
We have, for $1<\kappa<2$,
%
%
\begin{equation}
\Ep\bigl[\Var_\omega(\tau_{\mathrm{IA}})\bigr]=o_n\bigl
(n^{2/\kappa}
\bigr) \label{eqnlimevar},
\end{equation}
and, for $1\leq\kappa<2$,
\[
\Ep\bigl[\Var_\omega\bigl(\tau_{\mathrm{IA}}'\bigr)
\bigr]=o_n\bigl(n^{2/\kappa}\bigr).
\]
\end{lemma}

We recall that the second bound is only introduced to settle the case
$\kappa=1$; it would suffice for $1<\kappa<2$ as well, but introduces
unnecessary complication. The computations being very close for $\tau
_{\mathrm{IA}}$
and $\tau_{\mathrm{IA}}'$, we will write below the proof for $\tau
_{\mathrm{IA}}$ and indicate line
by line where changes happen for $\tau_{\mathrm{IA}}'$. Let us stress
that, when
dealing with $\tau_{\mathrm{IA}}'$, all the indicator functions
$\indic{H_\cdot<
h_n}$ (which define the small valleys) would refer to the original
heights, while all the potentials $V(\cdot)$ appearing along the
computation (which come from quenched expectations of times spent by
the walk) would refer to the modified environment.

\begin{pf*}{Proof of Lemma~\ref{lemvartia}}
We have
%
%
\begin{equation}
\label{eqntauiaexcursions} \tau_{\mathrm{IA}}=\sum_{p=0}^{n-1}
\tau(e_p,e_{p+1})\indic{H_p<
h_n},
\end{equation}
and by the Markov property, the above times are independent under
$P_{o,\omega}$. Hence,
\[
\Var_\omega(\tau_{\mathrm{IA}})=\sum_{p=0}^{n-1}
\Var_\omega\bigl(\tau(e_p,e_{p+1})\bigr)\indic
{H_p< h_n}.
\]
Under $\Pp$, the distribution of the environment seen from $e_p$ does
not depend on~$p$, hence,
%
%
\begin{equation}
\label{eqvar1} \Ep\bigl[\Var_\omega(\tau_{\mathrm{IA}})\bigr]=n
\Ep\bigl[
\Var_\omega\bigl(\tau(e_1)\bigr)\indic{H< h_n}
\bigr].
\end{equation}

We use formula \eqref{eqnalili2}:
%
%
\begin{equation}
\label{eqvarzeitouni1}
\Var_\omega\bigl(\tau(e_1)\bigr)
\indic{H< h_n}\leq16\sum_{z'\leq z\leq
x\leq
y\leq e_1, 0\leq y}
\ee^{{V}(y)+{V}(x)-{V}(z)-{V}(z')}\indic{H< h_n}.\hspace*{-40pt}
\end{equation}
Let us first consider the part of the sum where $x\geq0$. By noting
that the indices satisfy $z'\leq x$ and $z\leq y$, this part is seen to
be less than\break $(M'_1 M_2\ee^H)^2\indic{H<h_n}$. Lemma \ref
{lemrenewalestimates} shows that its expectation is smaller than $C\ee
^{(2-\kappa)h_n}$. \textit{For $\tau_{\mathrm{IA}}'$\textup{:} The same
holds\textup{,} because of the
remark on page \textup{\pageref{rem:e'r}}.}

It remains to deal with the indices $x<0$. This part rewrites as
%
%
\begin{equation}
\label{eqnjneg} \sum_{z',z\leq x< 0}\ee^{V(x)-V(z)-V(z')}\cdot
\sum_{0\leq y<e_1}\ee^{V(y)}\indic{H<h_n}.
\end{equation}
Since $V_{|\Z_+}$ and $V_{|\Z_-}$ are independent under $P$, so are the
two above factors. The second one equals $\ee^H M_2\indic{H<h_n}$. Let
us split the first one according to the excursion $[e_{p-1},e_p)$
containing $x$; it becomes
%
%
\begin{equation}
\label{eqn78648} \sum_{p\leq0}\ee^{-V(e_{p-1})}\sum
_{e_{p-1}\leq x<e_p}\ee^{V(x)-V(e_{p-1})} \biggl(\sum
_{z\leq x} \ee^{-(V(z)-V(e_{p-1}))} \biggr)^2.
\end{equation}
We have by definition $V(e_{p-1})\geq V(e_p)$ and, under $\Pp$,
$V(e_p)$ is independent of $(V(e_p+x)-V(e_{p}))_{x\leq0}$ and thus of
$(V(e_{p-1}+x)-V(e_{p-1}))_{x\leq e_{p}-e_{p-1}}$,\vspace*{1pt} which
has same
distribution as $(V(x))_{x\leq e_1}$. Therefore, the expectation of
\eqref{eqn78648} with respect to $\Pp$ is less than
\begin{eqnarray*}
&&\sum_{p\leq0} \Ep\bigl[\ee^{-V(e_p)}\bigr]\Ep
\biggl[\sum_{0\leq x<e_1}\ee^{V(x)} \biggl(\sum
_{z\leq x}\ee^{-V(z)} \biggr)^2
\biggr]
\\
&&\qquad\leq\bigl(1-E\bigl[\ee^{V(e_1)}\bigr]\bigr)^{-1}\Ep
\bigl[
\ee^H\bigl(M'_1\bigr)^2M_2
\bigr].
\end{eqnarray*}
Thus the expectation of \eqref{eqnjneg} with respect to $\Pp$ is
bounded by
\[
\bigl(1-E\bigl[\ee^{V(e_1)}\bigr]\bigr)^{-1}\Ep\bigl[
\ee^H\bigl(M'_1\bigr)^2M_2
\bigr]\Ep\bigl[\ee^H M_2\indic{H<h_n}\bigr].
\]
From Lemma~\ref{lemrenewalestimates}, we conclude that this term is
less than a constant if $\kappa>1$. The part corresponding to $x\geq0$
therefore dominates; this finishes the proof of \eqref{eqnlimevar}.
\textit{For $\tau_{\mathrm{IA}}'$\textup{:} The first factor is $(1-E[\ee
^{V(e_1)}|H<h_n])^{-1}$\textup{,} which is uniformly bounded because it
converges to $(1-E[\ee^{V(e_1)}])^{-1}<\infty$ and\textup{,} using Lemma \ref
{lemrenewalestimates}\textup{,} the two other factors are each bounded by a
constant if $\kappa>1$ and by $C h_n$ if $\kappa=1$\textup{;} cf. again
the remark page \textup{\pageref{rem:e'r}}. Thus\textup{,} the part corresponding to
$x\geq0$ still dominates in this case.}

We have proved $\Ep[\Var_\omega(\tau_{\mathrm{IA}})]\leq C n\ee
^{(2-\kappa)h_n}$.
Since $n\ee^{(2-\kappa)h_n}=\frac{n^{2/\kappa}}{(\log n)^{2-\kappa}}$,
this concludes the proof of \eqref{eqnlimevar}.\vspace*{-2pt}
\end{pf*}

\subsection{A subsequent lemma}

The proof of \eqref{eqnlimevar} entails the following bound for the crossing time of one low excursion.\vspace*{-2pt}
%
%
\begin{lemma}\label{lemtaue1}
For all $h>0$ we have, if $1<\kappa<2$,
\[
\Ep\bigl[E_\omega\bigl[\tau(e_1)^2\bigr]
\indic{H<h}\bigr]\leq C \ee^{(2-\kappa)h},
\]
and similarly for $(\Ep)'$ if $1\leq\kappa<2$.\vspace*{-2pt}
\end{lemma}

\begin{pf}
We have $E_\omega[\tau(e_1)^2]=\Var_\omega(\tau(e_1))+E_\omega
[\tau
(e_1)]^2$. Equation \eqref{eqvar1} and the remainder of the proof of
\eqref{eqnlimevar} give
\[
\Ep\bigl[\Var_\omega\bigl(\tau(e_1)\bigr)\indic{H< h}\bigr]
\leq C\ee^{(2-\kappa)h}.
\]
Together with Lemma~\ref{lemboundEe}, this concludes the proof.\vspace*{-2pt}
\end{pf}

\section{\texorpdfstring{Proof of Theorems \protect\ref{tquenched2} and \protect\ref{tquenched3}}
{Proof of Theorems 1 and 2}}
\label{secThm2-4}

Note that we first prove the results under $\Pp$. We will also
prove \eqref{tequenched} as a tool.\vspace*{-2pt}

\subsection{Joint coupling}

Extending what we did in Section~\ref{subseccoupling}, we introduce
an i.i.d. family $(\eg_i)_{i\geq1}$ of exponential random variables of
parameter $1$ and define, for $i\geq1$,
\[
N_i\defeq\biggl\lfloor-\frac{1}{\log(1-p_i(\omega))}\eg_i
\biggr\rfloor,
\]
where $p_i(\omega)=P_{b_i,\omega}(\tau(b_i)<\tau(d_i))$. Since, by
the Markov property, the numbers of returns to $b_i$ before the walk
reaches $d_i$ are independent given $\omega$, conditioning these
numbers to be equal to $N_i$ realizes a coupling, as in Section \ref
{subseccoupling}.\vspace*{-2pt}

\subsection{Reduction to one valley}
The above coupling enables us to give the following bound:
\begin{eqnarray*}
&&\W\Biggl(\tau(e_n)-E_\omega\bigl[\tau(e_n)
\bigr],\sum_{i=1}^{K_n}E_\omega
\bigl[\tau(b_i,d_i)\bigr]\egb_i \Biggr)
\\
&&\qquad\leq E_\omega\Biggl[\Biggl|\tau(e_n)-E_\omega\bigl[
\tau(e_n)\bigr]-\sum_{i=1}^{K_n}E_\omega
\bigl[\tau(b_i,d_i)\bigr]\egb_i \Biggr|
\Biggr]
\\
& &\qquad\leq E_\omega\bigl[ \bigl|\tau_{\mathrm{IA}}-E_\omega
[\tau_{\mathrm{IA}}] \bigr| \bigr]+\sum
_{i=1}^{K_n}E_\omega\bigl[ \bigl|
\tau(b_i,d_i)-E_\omega\bigl[\tau
(b_i,d_i)\bigr]\eg_i \bigr| \bigr],
\end{eqnarray*}
where $\tau_{\mathrm{IA}}$ is defined in Section \ref
{secinterarrival} (note that
for the $K_n$ high excursions the centerings simplify). We deduce, for
all $\delta\Sg0$,
\begin{eqnarray*}
&& \Pp\Biggl(\W\Biggl(\tau(e_n)-E_\omega\bigl[
\tau(e_n)\bigr],\sum_{i=1}^{K_n}E_\omega
\bigl[\tau(b_i,d_i)\bigr]\egb_i \Biggr)
\Sg\delta n^{1/\kappa} \Biggr)
\\[-2pt]
&&\qquad\leq\Pp\biggl(E_\omega\bigl[ \bigl|\tau_{\mathrm
{IA}}-E_\omega[\tau_{\mathrm{IA}}] \bigr| \bigr]
\Sg\frac
{\delta}{2}n^{1/\kappa} \biggr)
\\[-2pt]
& &\qquad\quad{}+ \Pp\Biggl(\bigcup_{p=0}^{n-1} \biggl\{
E_\omega\bigl[ \bigl|\tau(e_p,e_{p+1})-E_\omega
\bigl[\tau(e_p,e_{p+1})\bigr]\eg_p \bigr| \bigr]
\indic{H_p\geq h_n}\\[-2pt]
&&\hspace*{235pt}\qquad{}\geq\frac{\delta}{2K_n}n^{1/\kappa}
\biggr\} \Biggr).
\end{eqnarray*}
By Proposition~\ref{proptauia}, the first term is known to converge
to $0$ as $n\to\infty$ (using for $\kappa=1$ the same reduction as in
Section~\ref{secinterarrival}). By Lemma~\ref{lemenegatifs}, the
last term is bounded by
\begin{eqnarray*}
P \bigl(K_n\geq2(\log n)^\kappa\bigr)+n \Pp
\biggl(E_\omega\bigl[ \bigl|\tau-E_\omega[\tau]\eg\bigr| \bigr
]\geq
\frac{\delta}{4(\log n)^\kappa
}n^{1/\kappa},H\geq h_n \biggr),
\end{eqnarray*}
where $\tau$ and $\eg$ stand for $\tau(e_1)$ ($=\tau(b_1,d_1)$ on
$\{
H\geq h_n\}$) and $\eg_1$. By the proof of Proposition \ref
{propoverlap}, the first probability goes to $0,$ when $n$ tends to
infinity. As for the other probability, it follows from
Proposition \ref
{propexitonevalley} with $t=n^{1/\kappa}$ that it is $o(n^{-1})$.

This yields, under $\Pp$,
%
%
\begin{equation}\label{eqthmquencheden}
\W\Biggl(\frac{\tau(e_n)-E_\omega[\tau(e_n)]}{n^{1/\kappa
}},\frac
{1}{n^{1/\kappa}}\sum
_{i=1}^{K_n}E_\omega\bigl[
\tau(b_i,d_i)\bigr]\egb_i \Biggr)
\limites{(p)} {n}0,
\end{equation}
which is the statement of \eqref{tequenched} along the random
subsequence $x=e_n$, and under $\Pp$ instead of $P$.
Before proceeding to the interpolation from $e_n$ to any $x$, let us
show how the statements of Theorems~\ref{tquenched2} and \ref
{tquenched3} can be quickly deduced from \eqref{eqthmquencheden},
modulo the same restriction.

\subsection{Addition of small excursions and independence of the high
ones}\label{subsecajoutpetites}

More specifically, if (with a convenient abuse of notation) we extend
the i.i.d. sequence $(\egb_i)_{i\geq1}$ to an i.i.d. sequence $(\egb
_p)_{p\geq0}$ such that $\egb_i=\egb_p$ for $p=\sigma(i)$, the only
addition in Theorem~\ref{tquenched2} is the following term which we
shall prove is negligible:
%
%
\begin{equation}
\frac{1}{n^{1/\kappa}}\sum_{p=0}^{n-1}Z_p
\indic{H_p\Sl h_n}\egb_p,
\label{eqpetitescentreesponderees}\vadjust{\goodbreak}
\end{equation}
where we define
\[
Z_i\defeq E_\omega\bigl[\tau(e_i,e_{i+1})
\bigr],\qquad i \ge0.
\]
Note that $(Z_i)_{i\geq0}$ is a stationary sequence under $\Pp$;
cf. Lemma~\ref{lemenegatifs}.

For $0\Sl\kappa\Sl1$, it suffices to note that the $L^1(P_\omega)$-norm
of this term is bounded by
$n^{-1/\kappa}E_\omega[\tau_{\mathrm{IA}}]$
(since $E_\omega[|\egb_p|]=2/\ee\Sl1$), which converges to 0 in
$L^1(P)$ and thus in probability in this case; cf. after
Proposition~\ref{proptauia}.

For $1\Sl\kappa\Sl2$, let us write that the $L^1(P_\omega)$-norm
of \eqref{eqpetitescentreesponderees} is bounded, using
the Cauchy--Schwarz inequality, by
\[
\frac{1}{n^{1/\kappa}}\Var_\omega\Biggl(\sum_{p=0}^{n-1}
Z_p \indic{H_p\Sl h_n}\egb_p
\Biggr)^{1/2}=\frac{1}{n^{1/\kappa}} \Biggl(\sum
_{p=0}^{n-1}Z_p^2
\indic{H_p\Sl h_n} \Biggr)^{1/2},
\]
hence,
\begin{eqnarray*}
\Pp\Biggl(E_\omega\Biggl[ \Biggl|\frac{1}{n^{1/\kappa}}\sum
_{p=0}^{n-1}Z_p \indic{H_p
\Sl h_n}\egb_p \Biggr| \Biggr]\geq\delta\Biggr)\leq
\frac
{1}{\delta^2 n^{2/\kappa}}n\Ep\bigl[E_\omega[\tau]^2\indic{H\Sl
h_n} \bigr].
\end{eqnarray*}
Lemma~\ref{lemboundEe} shows that the last expectation is less than
$C n^{\frac{2}{\kappa}-1}(\log n)^{-(2-\kappa)}$ so that the right-hand
side converges to $0$.

For $\kappa=1$, we do the same as for $\kappa\Sg1$, by means of the
reduction to the modified environment (cf. Section \ref
{subsecreducsmall}): the decomposition $\tau_{\mathrm{IA}}=\tau
_{\mathrm{IA}}'-\widetilde{\tau}_{\mathrm{IA}
}+\widetilde{\tau}_{\mathrm{IA}}'$ of~\eqref{eqtiaprime} induces a
decomposition
similar to \eqref{eqpetitescentreesponderees} (with the only addition
of quenched expectations and weights). The terms corresponding to
$\widetilde{\tau}_{\mathrm{IA}}$ and $\widetilde{\tau}_{\mathrm{IA}}'$ are neglected using their
first moment by the results \eqref{eqnlimtiatilde} and \eqref
{eqnlimtiatildee} in Section~\ref{subsecreducsmall}, thus
reducing the problem to the modified environment, where Lemma \ref
{lemtaue1} applies. This would conclude the proof of Theorem \ref
{tquenched2}, up
to the previous restrictions.

To deduce Theorem~\ref{tquenched3} (along the random subsequence
$x=e_n$ and under $\Pp$) from (\ref{eqthmquencheden}), we have to
replace $Z_{\sigma(i)}=E_\omega[\tau(b_i,d_i)]$, $i=1,\ldots,K_n$, by
independent terms having the same distribution, and to add new terms
corresponding to small excursions, just like above but independent of
each other. Note that the new independent terms $\widehat{Z}_p$ will
depend on $n$, even though their distribution does not, which explains
the wording of Theorem~\ref{tquenched3}.

To this aim, let us enlarge the probability space $(\Omega\times\Z
^\N
,\mathcal{B},\bPp)$ in order to introduce a sequence $(\omega
^{(p)},(X^{(p)}_t)_{t\in\N})_{p\geq0}$ of environments and random
walks coupled with $\omega$ in the following way, for $p\geq0$:
\begin{longlist}[(1)]
\item[(1)] if $H_p<h_n$, then $\omega^{(p)}$ is an independent
environment sampled according to the distribution $\Pp(\cdot|H<h_n)$;
\item[(2)] if $H_p\geq h_n$, that is, $p=\sigma(i)$ for some $i\geq1$,
then $\omega^{(p)}$ is built from the piece of $\omega$ from
$d_{i-1}+1$ to $d_i$, translated so that $b_i$ is now at 0, and
bordered by independent environments with law $P$ on the right and law
$\Pp(\cdot|H_{-1}\geq h_n,V_{|\Z_-}\geq-A_i)$ on the left where
$A_i\defeq V(d_{i-1})-V(b_i)$ (function of $\omega$);
\item[(3)] for all $p\geq0$, conditionally on $\omega^{(p)}$,
$(X^{(p)}_t)_{t\in\N}$ has law $P_{\omega^{(p)}}$.\vadjust{\goodbreak}
\end{longlist}
Due to the independence between the excursions of $\omega$ under $\Pp$,
the sequence $(\omega^{(p)})_{p\geq0}$ is seen to be independent.
Furthermore, for every $p\geq0$, the construction ensures that $\omega
^{(p)}$ follows the law $\Pp$. We will denote with a superscript
${}^{(p)}$ the quantities relative to $\omega^{(p)}$ instead of
$\omega$.

We may thus introduce
\[
\widehat{Z}_p\defeq E_\omega\bigl[ \tau^{(p)}
\bigl(e_1^{(p)}\bigr)\bigr],\qquad p\geq0,
\]
which is defined as $Z_1(\defeq E_\omega[ \tau(e_1)]) $ but relative
to $(\omega^{(p)},X^{(p)})$ instead of $(\omega,X)$. By the previous
claims, $(\widehat{Z}_p)_{p\geq0}$ is a sequence of i.i.d. random
variables distributed as $Z_1$ under $\Pp$.

For $i\geq1$, to compare $Z_{\sigma(i)}$ with $\widehat{Z}_{\sigma
(i)}$, we further decompose
\[
Z_{\sigma(i)}\defqe\widetilde{Z}_{\sigma(i)}+Z^*_{\sigma(i)}
\quad\mbox{and}\quad\widehat{Z}_{\sigma(i)}\defqe\widehat
{\widetilde
{Z}}_{\sigma
(i)}+
\widehat{Z}^*_{\sigma(i)},
\]
where we let $\widetilde{Z}_{\sigma(i)}\defeq E_\omega[\tautx
{a_i}(b_i,d_i)]$ and similarly, $\widehat{\widetilde{Z}}_{\sigma(i)}$
is defined as\break $E_\omega[\tautx{e_{-D_n}}(0,e_1)]$ with respect to
$\omega^{(\sigma(i))}$ instead of $\omega$, so that $Z^*_{\sigma(i)}$
is the quenched expectation of the time to go from $b_i$ to $d_i$ for a
random walk reflected at $a_i$ and thus only depends on the environment
between $a_i$ and $d_i$. Using this last remark, it is important to
note that, on the event $\operatorname{NO}(n)$ (cf. Proposition~\ref
{propoverlap}),
$Z^*_{\sigma(i)}$ and $\widehat{Z}^*_{\sigma(i)}$ are equal for
$i=1,\ldots,K_n$. Indeed, since $P(\operatorname{NO}(n))\to_n1$,
this gives us directly
%
%
\begin{equation}
\label{eqthm2laststep} \W\Biggl(\frac{1}{n^{1/\kappa}} \sum
_{i=1}^{K_n} Z_{\sigma(i)}^* \egb_{\sigma(i)},
\frac{1}{n^{1/\kappa}}\sum_{i=1}^{K_n} \widehat
{Z}_{\sigma
(i)}^* \egb_{\sigma(i)} \Biggr)\limites{(p)} {n}0.
\end{equation}
In addition, Proposition~\ref{proplimtaut} and the triangular
inequality give
%
%
\begin{equation}
\W\Biggl(\frac{1}{n^{1/\kappa}} \sum_{i=1}^{K_n}
\widetilde{Z}_{\sigma
(i)} \egb_{\sigma(i)},0 \Biggr)\limites{(p)} {n}0.
\end{equation}
It remains to prove that the same holds for $\widehat{\widetilde
{Z}}_{\sigma(i)}$ in order to get \eqref{eqthmquencheden} with
$\widehat{Z}^{(\sigma(i))}$ in place of $Z^{(\sigma(i))}$. And finally,
Theorem~\ref{tquenched3} will be proved (under the above-mentioned
restrictions) if, furthermore, the small independent excursions may be
harmlessly introduced, that is, if
%
%
\begin{equation}
\W\Biggl(\frac{1}{n^{1/\kappa}} \sum_{p=0}^{n-1}
\widehat{Z}_p \egb_p\indic{H^{(p)}<h_n},0
\Biggr)\limites{(p)} {n}0.
\end{equation}
These two facts are given by the following lemma.

%
\begin{lemma}\label{lemzhat}
We have, under $\Pp,$
\[
\frac{1}{n^{1/\kappa}}\sum_{p=0}^{n-1}
\widehat{\widetilde{Z}}_p \indic{H^{(p)}\geq
h_n} \limites{(p)} {n} 0\label{eqnlimtauthat}
\]
and
%
%
\begin{equation}
\frac{1}{n^{1/\kappa}} \sum_{p=0}^{n-1}
\widehat{Z}_p \indic{H^{(p)}<h_n}
\egb_p \limites{(p)} {n}0.\label{eqntauiahat}
\end{equation}
\end{lemma}

\begin{pf}
These results follow, respectively, from the proofs of Proposition \ref
{proplimtaut} and (\ref{eqpetitescentreesponderees}), made easier by
the independence of the random variables $\widehat{Z}_0,\ldots
,\widehat
{Z}_{n-1}$.
More precisely, the proof of Proposition~\ref{proplimtaut} holds in
this i.i.d. context almost without a change, while the above derivation
of \eqref{eqpetitescentreesponderees} did not involve the correlation
between $Z_0,\ldots,Z_{n-1}$ in any way, hence, the proof may as well
be conducted for independent copies.
\end{pf}

\subsection{\texorpdfstring{Interpolation from $\tau(e_n)$ to $\tau(x)$}{Interpolation from tau(e_n) to tau(x)}}
We now replace
the subsequence $\tau(e_n)$ by the whole sequence $\tau(x)$. We write
the proof in the setting of Theorem~\ref{tquenched2}, from which the
other cases follow, up to very minor modifications.

Choose $\frac{1}{2}<\alpha<\min\{1,\frac{1}{\kappa}\}$. For $x\in
\N$,
we define the following event about the environment:
%
%
\begin{equation}
A_x\defeq\{e_{ \lfloor{(x-x^\alpha)}/{(E[e_1])}
\rfloor
}<x<e_{ \lfloor{(x+x^\alpha)}/{(E[e_1])} \rfloor} \}.
\end{equation}
Since $\alpha>\frac{1}{2}$, it follows from the central limit theorem,
applied to the i.i.d. sequence $(e_{n+1}-e_n)_n$, that
%
%
\begin{eqnarray}
P(A_x)\to1,\qquad x \to\infty. \label{eqnaxandbx}
\end{eqnarray}
Starting from the version of Theorem~\ref{tquenched2} we have obtained
so far, that is, for every $\delta>0$,
\[
\Pp\Biggl( \Biggl|\tau(e_n)-E_\omega\bigl[\tau(e_n)
\bigr]-\sum_{p=0}^{n-1} E_\omega
\bigl[\tau(e_p,e_{p+1})\bigr]\egb_p \Biggr|>\delta
n^{1/\kappa} \Biggr)\limites{} {n}0,
\]
the limit still holds along the \textit{deterministic} subsequences
\[
n^-_x\defeq\biggl\lfloor\frac{x-x^\alpha}{E[e_1]} \biggr\rfloor
\quad\mbox{and}\quad n^+_x\defeq\biggl\lfloor\frac{x+x^\alpha
}{E[e_1]} \biggr\rfloor,
\]
and according to \eqref{eqnaxandbx} it is legitimate to restrict
to the event $A_x$ in the above probability for $n=n_x^\pm$. From that
remark and $n^\pm_x\sim_x \frac{x}{E[e_1]}$, we conclude that the
result of Theorem~\ref{tquenched2} will follow from (under $\Pp$)
\[
\frac{1}{x^{1/\kappa}}E_\omega\bigl[ \bigl|\tau(x)-\tau
(e_{n^+_x})\bigr|
\bigr]\limites{(p)} {x}0,\qquad\frac{1}{x^{1/\kappa}} \bigl
|E_\omega
\bigl[\tau(x)
\bigr]-E_\omega\bigl[\tau(e_{n^+_x})\bigr]\bigr|\limites{(p)} {x}0,
\]
the corresponding limits for $n^-_x$ and
\[
\frac{1}{x^{1/\kappa}}\sum_{n^-_x\leq p\leq n^+_x}E_\omega\bigl[
\tau(e_p,e_{p+1})\bigr]\egb_p\limites{(p)}
{x}0.
\]
Of course, the second limit will follow from the first one.
Furthermore, on $A_x$ we have
\[
E_\omega\bigl[ \bigl|\tau(x)-\tau(e_{n^\pm_x}) \bigr| \bigr
]\leq
E_\omega\bigl[\tau(e_{n^+_x})-\tau(e_{n^-_x}) \bigr]=
\sum_{n^-_x\leq p<
n^+_x}E_\omega\bigl[
\tau(e_p,e_{p+1})\bigr],
\]
so that the three limits will come as a consequence of the following
application of the Markov inequality:
\begin{eqnarray*}
&&\Pp\biggl(\sum_{n^-_x\leq p\leq n^+_x}E_\omega\bigl[\tau
(e_p,e_{p+1})\bigr]>\delta x^{1/\kappa} \biggr)
\\
&&\qquad\leq P\bigl(\exists n^-_x\leq p\leq n^+_x,
H_p\geq h_x\bigr)+\frac
{n^+_x-n^-_x+1}{\delta x^{1/\kappa}}\Ep
\bigl[E_\omega\bigl[\tau(e_1)\bigr],H<h_x \bigr]
\\
&&\qquad\leq\frac{2x^\alpha}{E[e_1]}P(H\geq h_x)+\frac{2x^\alpha
+1}{\delta
x^{1/\kappa}}\Ep
\bigl[2M'_1M_2\ee^H,H<h_x
\bigr].
\end{eqnarray*}
By \eqref{iglehartthm} and $\alpha<1$, the first term goes to 0. By
Lemma~\ref{lemrenewalestimates} and since $\alpha<\frac1\kappa$, the
second term goes to 0 as well. This proves Theorem~\ref{tquenched2},
under $\Pp$.

\subsection{Conclusion}\label{subsecccl} Let us finally discuss the
change of probability from $\Pp$ to~$P$. In fact, it suffices to note
that the quenched expectation of the time spent on $\Z_-$ is finite
a.s. under $P$ and $\Pp$, which follows from \eqref{eqnzeitounip3}
and \eqref{eqnzeitounie3} (and $E[\log\rho]<0$) since this
expectation\vspace*{1pt} is seen to be equal to $E_{0,\omega}[\tau
(1)]P_{1,\omega
}(\tau(0)=\infty)^{-1}$. This ends the proof of \eqref{tequenched} and
Theorems~\ref{tquenched2} and~\ref{tquenched3}.

Note that the tail estimate \eqref{eqnqueuezh} of $\widehat{Z}_i$
(i.e., of $E_\omega[\tau(e_1)]$ under $\Pp$) given in Theorem \ref
{tquenched3}, while not being exactly a consequence of Lemma \ref
{lemtailz}, follows simply from it. Indeed,\vspace*{1pt} the
expression $E_\omega
[\tau(e_1)]=E_\omega[N]E_\omega[F]+E_\omega[G]=2\ee^H \widehat
{M}_1M_2+E_\omega[G]$, together with \eqref{eqboundM1} and
Lemma \ref
{lemmafsuccessbound}, gives the following lower and upper bounds, for
some $\alpha<1$ and $\delta>0$:
\begin{eqnarray*}
\Pp\biggl(2Z\geq\frac{t}{1+o(t^{-\delta})} \biggr) & \leq&\Pp
\bigl(E_\omega
\bigl[\tau(e_1)\bigr]\geq t\bigr)
\\
& \leq&\Pp\bigl(\Omega_t^c\bigr)+\Pp\biggl(2Z\geq
\frac{t-C(\log t)^4 t^\alpha
}{1+o(t^{-\delta})} \biggr),
\end{eqnarray*}
and $\Pp(\Omega_t^c)=o(t^{-\kappa})$ by Lemma~\ref{lemgoodenv}, hence,
with Lemma~\ref{lemtailz},
%
%
\begin{equation}
\label{eqqueueetau} \Pp\bigl(E_\omega\bigl[\tau(e_1)\bigr]
\geq t \bigr)\sim2^\kappa C_U t^{-\kappa
},\qquad t\to\infty.
\end{equation}

\section{\texorpdfstring{Proof of Corollary \protect\ref{coroll}}{Proof of Corollary 1}}
\label{secproofcoro}

We show here how Corollary~\ref{coroll} follows from Theorem \ref
{tquenched3}. With the notation of this theorem, it suffices to prove
\[
\mathscr{L} \Biggl(\frac{1}{x^{1/\kappa}}\sum_{p=1}^{x}
\widehat{Z}_p\egb_p\Big\rrvert(\widehat{Z}_p)_{p\geq1}
\Biggr)\limites{\mathrm{W}^1} {x} \mathscr{L} \Biggl(\sum
_{p=1}^\infty\xi_p \egb_p
\Big\rrvert(\xi_p)_{p\geq1} \Biggr)\qquad\mbox{in law,}
\]
where $(\widehat{Z}_p)_{p\geq1}$ are i.i.d., independent of $(\egb
_p)_{p\geq1}$, such that $P(\widehat{Z}_1>t)\sim2^{\kappa} C_U
t^{-\kappa}$,
and $(\xi_p)_{p\geq1}$ is a Poisson point process of intensity
$2^{\kappa} C_U \kappa u^{-(\kappa+1)} \,\d u $, independent of
$(\egb_p)_{p\ge1}$.
This reduction comes from the following easy property.

%
\begin{lemma}\label{lemprobaloi}
If random variables $(X_n)_n$, $(Y_n)_n$ and $Y$ take values in a
metric space $(E,d)$, $d(X_n,Y_n)\to_n 0$ in probability and $Y_n\to_n
Y$ in law imply $X_n\to_n Y$ in law.
\end{lemma}

Let us recall a simple result about order statistics of heavy-tailed
random variables.
%
%
\begin{proposition}
Let $(Z_i)_{i\geq1}$ be i.i.d. copies of a random variable $Z\geq0$
such that
%
%
\begin{equation}
\label{eqnqueuezlem} P(Z>t)\sim C_Zt^{-\kappa},\qquad t\to\infty,
\end{equation}
for some constant $C_Z>0$. For all $n\geq1$, denote by $Z^{(1)}_n\geq
\cdots\geq Z^{(n)}_n$ an ordering of the finite subsequence
$(Z_1,\ldots
,Z_n)$. Then we have, for every \mbox{$k\geq1$},
\[
\frac{1}{n^{1/\kappa}}\bigl(Z^{(1)}_n,\ldots,Z^{(k)}_n
\bigr)\limites{\mathrm{law}} {n} \bigl(\xi^{(1)},\ldots,\xi^{(k)}
\bigr),
\]
where $\xi^{(k)}=C_Z^{1/\kappa}(\fg_1+\cdots+\fg_k)^{-1/\kappa}$ for
$k\geq1$, $(\fg_k)_k$ being i.i.d. exponential random variables of
parameter $1$; cf. \eqref{eqndescrxi}.
\end{proposition}

\begin{pf}
Let $Y^{(i)}_n\defeq n C_Z (Z^{(i)}_n)^{-\kappa}$, and $Y_n=n
C_Z(Z_1)^{-\kappa}$. From \eqref{eqnqueuezlem} we deduce
$nP(Y_n\in
[a,b])\to_n b-a$ for all $0<a<b$. Then, for all $t_1,\ldots,t_k>0$,
\begin{eqnarray*}
&&P\bigl(t_1<Y^{(1)}_n<t_2<
Y^{(2)}_n<\cdots<t_k<Y^{(k)}_n
\bigr)\\
&&\qquad= n(n-1)\cdots\bigl(n-(k-1)+1\bigr)
P\bigl(Y_n\in[t_1,t_2]\bigr)\cdots
P\bigl(Y_n\in[t_{k-1},t_k]\bigr)\\
&&\qquad\quad{}\times P
\bigl(Y_n\notin[0,t_k]\bigr)^{n-k}
\\
&&\qquad\to_n (t_2-t_1)\cdots(t_k-t_{k-1})
\ee^{-t_k}
\\
&&\qquad= P(t_1<\fg_1<t_2<\fg_1+
\fg_2<\cdots<t_k<\fg_1+\cdots+
\fg_k),
\end{eqnarray*}
by a simple computation, from where the proposition follows.
\end{pf}

\newcommand{\Zt}{\widetilde{Z}}
\newcommand{\Zh}{\widehat{Z}}
\newcommand{\xit}{\widetilde{\xi}}

Thanks to the previous lemma and Skorohod's representation theorem,
there exists a copy $(\xit^{(p)})_{p\geq1}$ of $(\xi^{(p)})_{p\geq1}$
and, for all $k\geq1$, there exist random variables $(\Zt
^{(1)}_{k,n},\ldots,\Zt^{(k)}_{k,n})_{n\geq k}$
such that (borrowing notation from the lemma) for every $n\geq k$ $(\Zt
^{(1)}_{k,n},\ldots,\Zt^{(k)}_{k,n})$ is a copy of $(\Zh
^{(1)}_n,\ldots
,\Zh^{(k)}_n)$ and
\[
\frac{1}{n^{1/\kappa}}\bigl(\Zt^{(1)}_{k,n},\ldots,
\Zt^{(k)}_{k,n}\bigr)\limites{(p)} {n} \bigl(
\xit^{(1)},\ldots,\xit^{(k)}\bigr).
\]
We chose $(\xit^{(p)})_{p\geq1}$ not depending\vspace*{-1.5pt} on $k$ to ease notation
but this is unessential since we only need to understand the
convergences in probability $X_n{\displaystyle\mathop{\longrightarrow
}^{(p)}_{n}} X$ as properties of the
law of $(X_n,X)$ for every $n$, no matter on which space $\Omega_n$
this couple is defined.

We may also introduce additional random variables $(\Zt
^{(k+1)}_{k,n},\ldots,\Zt^{(n)}_{k,n})_{n\geq1}$
such that for every $n$ $(\Zt^{(1)}_{k,n},\ldots,\Zt^{(n)}_{k,n})$
is a
copy of $(\Zh^{(1)}_n,\ldots,\Zh^{(n)}_n)$.

Then, by a diagonal argument, we can define $(\Zt^{(p)}_n)_{1\leq
p\leq n}$
such that, for every $n$, $(\Zt^{(p)}_n)_{1\leq p\leq n}$ is a copy of
$(Z^{(1)}_n,\ldots,Z^{(n)}_n)$ and, \textit{for every $k$},
%
%
\begin{equation}
\label{eqncvps} \frac{1}{n^{1/\kappa}}\bigl(\Zt^{(1)}_{n},
\ldots,\Zt^{(k)}_{n}\bigr)\limites{(p)} {n} \bigl(
\xit^{(1)},\ldots,\xit^{(k)}\bigr).
\end{equation}
Indeed,
there is an increasing sequence $(N(k))_k$ such that for all $k\ge1$,
for $n\geq N(k)$,
\[
P \biggl( \biggl\|\frac{1}{n^{1/\kappa}}\bigl(\Zt
^{(1)}_{k,n},\ldots,
\Zt^{(k)}_{k,n}\bigr)-\bigl(\xit^{(1)},\ldots,
\xit^{(k)}\bigr) \biggr\|_1>\frac1k \biggr)< \frac1k
\]
(hence, the same bound also holds for the first $k'\leq k$ components)
and then we define, for $n\geq N(1)$ and $1\leq p\leq n$, $\Zt
^{(p)}_n=\Zt^{(p)}_{k,n},$ where $k$ is given by $N(k)\leq n<N(k+1)$;
and, for instance, $\Zt^{(p)}_n=\Zt^{(p)}_{1,n}$ when $1\leq p\leq
n<N(1)$. This is easily seen to satisfy \eqref{eqncvps}.

We have, for all $n\geq k$,
%
%
\begin{eqnarray}\label{eqndecoupw1}
&& W^1_{\Zt,\xit} \Biggl(\sum_{p=1}^n
\frac{\Zt^{(p)}_n}{n^{1/\kappa
}}\egb_p,\sum_{p=1}^\infty
\xit^{(p)}\egb_p \Biggr)\nonumber\\
&&\qquad\leq E_{\Zt,\xit} \Biggl[
\Biggl\llvert\sum_{p=1}^n
\frac{\Zt^{(p)}_n}{n^{1/\kappa}}\egb_p-\sum_{p=1}^\infty
\xit^{(p)}\egb_p\Biggr\rrvert\Biggr]
\nonumber
\\
&&\qquad\leq E_{\Zt} \Biggl[\Biggl\llvert\sum_{p=k+1}^{n}
\frac
{\Zt^{(p)}_n}{n^{1/\kappa}}\egb_p\Biggr\rrvert\Biggr] + E_{\Zt
,\xit}
\Biggl[\Biggl\llvert\sum_{p=1}^k \biggl(
\frac{\Zt^{(p)}_n}{n^{1/\kappa}}-\xit^{(p)} \biggr)\egb_p\Biggr
\rrvert
\Biggr] \\
&&\qquad\quad{}+ E_{\xit} \Biggl[\Biggl\llvert\sum
_{p=k+1}^\infty\xit^{(p)}\egb_p
\Biggr\rrvert\Biggr]
\nonumber\\
&&\qquad\leq\sqrt{\sum_{p=k+1}^{n}
\biggl(\frac{\Zt^{(p)}_n}{n^{1/\kappa}} \biggr)^2} + \sum
_{p=1}^k\biggl\llvert\frac{\Zt^{(p)}_n}{n^{1/\kappa}}-
\xit^{(p)}\biggr\rrvert+ \sqrt{\sum
_{p=k+1}^\infty\bigl(\xit^{(p)}
\bigr)^2},\nonumber
\end{eqnarray}
using $E[|\egb_p|]=2/ \ee\le1$ and the inequality $E[|W|]^2\leq
E[W^2]=\operatorname{ Var}(W)$ for any centered random variable $W$.
Let $\eps_k\Sg0$ be such that $k^{-1/\kappa}\ll\eps_k\ll1$,\vfill\eject when
$k\to
\infty.$
Since $\Zh^{(k)}_n\geq\Zh^{(p)}_n$ for $p\geq k$,
\begin{eqnarray*}
&& P \Biggl(\sqrt{\sum_{p=k+1}^{n}
\biggl(\frac{\Zh^{(p)}_n}{n^{1/\kappa}} \biggr)^2}\geq\frac
\delta3 \Biggr)
\\
&&\qquad\leq P \biggl(\frac{\Zh^{(k)}_n}{n^{1/\kappa}}\geq\eps_k
\biggr
)+P \Biggl(\sum
_{p=1}^n \biggl(\frac{\Zh_p}{n^{1/\kappa}}
\biggr)^2\ind_{\{
{\Zh
_p}/{n^{1/\kappa}}\Sl\eps_k\}}\geq\biggl(\frac{\delta}{3}
\biggr)^2 \Biggr)
\\
&&\qquad\leq P \biggl(\frac{\Zh^{(k)}_n}{n^{1/\kappa}}\geq\eps_k
\biggr)+
\frac
{9}{\delta^2}n E \biggl[ \biggl(\frac{\Zh_1}{n^{1/\kappa}} \biggr)^2
\ind_{\{
{\Zh_1}/{n^{1/\kappa}}\Sl\eps_k\}} \biggr],
\end{eqnarray*}
hence, using \eqref{eqncvps} and \eqref{eqnqueuezh},
for all $\delta>0$,
\[
\mathop{\limsup}_n P \Biggl(\sqrt{\sum
_{p=k+1}^{n} \biggl(\frac{\Zh^{(p)}_n}{n^{1/\kappa}}
\biggr)^2}\geq\frac\delta3 \Biggr) \leq P\bigl(\xi^{(k)}
\geq\eps_k\bigr)+\frac{9}{\delta^2}\frac
{2C}{2-\kappa}
\eps_k^{1-\kappa/2}\defqe\varphi_\delta(k),
\]
where $C>C_Z\defeq2^\kappa C_U$ is arbitrary. Note that $\varphi
_\delta
(k)\to_k 0$ due to the choice of $\eps_k$ and to \eqref{eqnequivxi}.
We also have, respectively, because of \eqref{eqncvps} and of $\sum_p
(\xi^{(p)})^2\Sl\infty$ a.s. [cf. \eqref{eqnequivxi}],
\begin{eqnarray*}
P \Biggl(\sum_{p=1}^k\biggl\llvert
\frac{\Zt^{(p)}_n}{n^{1/\kappa}}-\xit^{(p)}\biggr\rrvert\geq
\frac{\delta}{3} \Biggr)
\limites{} {n}0 \quad\mbox{and}\quad P \Biggl(\sqrt{\sum
_{p=k+1}^\infty\bigl(\xi^{(p)}
\bigr)^2}\geq\frac
{\delta}{3} \Biggr)\limites{} {k}0.
\end{eqnarray*}
Denote by $\psi_\delta(k)$ the latter probability. Thus, from \eqref
{eqndecoupw1}, for all $\delta>0$,
\begin{eqnarray*}
\mathop{\limsup}_n P \Biggl(W^1_{\Zt,\xit} \Biggl(\sum
_{p=1}^n\frac{\Zt^{(p)}_n}{n^{1/\kappa}}
\egb_p,\sum_{p=1}^\infty
\xit^{(p)}\egb_p \Biggr)\geq\delta\Biggr)\leq
\varphi_\delta(k)+\psi_\delta(k) \to_k 0.
\end{eqnarray*}
Thanks to our diagonal argument, the left-hand side does not depend on
$k$. Thus,
\[
\mathscr{L} \Biggl(\sum_{p=1}^n
\frac{\Zt^{(p)}_n}{n^{1/\kappa
}}\egb_p \Big|\bigl(\Zt_n^{(p)}
\bigr)_{1\leq p\leq n} \Biggr)\limites{W^1} {n}\mathscr{L} \Biggl(
\sum_{p=1}^\infty\xit^{(p)}
\egb_p \Big|\xit\Biggr) \qquad\mbox{in probability},
\]
and therefore in law. Since the convergence in law only deals with the
laws of $\Zt_n$ for $n\geq1$ and of $\xit$ (and not on their
coupling), this concludes the proof of Corollary~\ref{coroll}.

Finally, we mention that the expression of the parameter $\lambda$
obtained for Dirichlet environments [i.e., when $\omega_0$ follows a
distribution $\operatorname{ Beta}(\alpha,\beta)$ with $0\Sl\alpha
-\beta\Sl2$]
can be easily deduced from a computation of $C_K$ by Chamayou and Letac
\cite{chamayou-letac} (see~\cite{limitlaws} for more details).

\begin{appendix}
\section*{Appendix}\label{app}

\subsection{\texorpdfstring{Proofs of Lemmas \protect\ref{lemtailz} and \protect\ref{lemrenewalestimates}}
{Proofs of Lemmas 2 and 3}}
\label{subsecprooflemma-renewalestimate}
\mbox{}
\begin{pf*}{Proof of Lemma~\ref{lemtailz}}
Compared to \eqref{eqnqueuezcondition}, it appears sufficient to
prove that $\Pp(Z>t,S>H)=o(t^{-\kappa})$, which is understood as
follows: when $Z$ is large, the height $H$ of the first excursion tends
to be large as well, while the other excursions are independent of $Z$,
hence, $H$ is likely to be the maximum $S$ of $V$ over all of $\Z_+$.
More precisely: first, for $\ell_t>0$,
\[
\Pp(Z>t, H<\ell_t)\leq\Pp\bigl(M_1 M_2>t
\ee^{-\ell_t}\bigr)\leq\frac{\Ep
[(M_1M_2)^2]}{(t\ee^{-\ell_t})^2},
\]
and all moments of $M_1M_2$ are finite under $\Pp$ [indeed we have
$M_2\leq e_1$, $M_1\leq e_1+R_-$ and the random variables $e_1$ and
$R_-$ have all moments finite under $\Pp$; cf. Section \ref
{subsecrenewal} and \eqref{eqnr-}]. Thus, if (recalling that $\kappa
<2$) we choose $\ell_t$ such that $\ell_t\to\infty$ and $t^\kappa
=o(t^2\ee^{-2\ell_t})$ as $t\to\infty$, we have $\Pp(Z>t,H<\ell
_t)=o(t^{-\kappa})$. On the other hand, $Z$ is independent of
$S'\defeq
\sup_{x\geq e_1}(V(x)-V(e_1))$ which is larger than $S$ on the event
$\{
H<S\}$, hence,
\begin{eqnarray*}
\Pp(Z>t, H<S) & =& \Pp(Z>t, H\geq\ell_t, H<S) + o
\bigl(t^{-\kappa}\bigr)
\\
& \leq&\Pp(Z>t)\Pp\bigl(S'>\ell_t\bigr)+o
\bigl(t^{-\kappa}\bigr)
\\
& =&\Pp(Z>t)o(1)+o\bigl(t^{-\kappa}\bigr),
\end{eqnarray*}
as $t\to\infty$, such that, using \eqref{eqnqueuezcondition},
\begin{eqnarray*}
\Pp(Z>t) & =& \Pp(Z>t,H=S)+\Pp(Z>t,H<S)
\\
& = &C_U t^{-\kappa}+o\bigl(t^{-\kappa}\bigr)+
\Pp(Z>t)o(1)+o\bigl(t^{-\kappa}\bigr),
\end{eqnarray*}
which implies the lemma.
\end{pf*}

\begin{pf*}{Proof of Lemma~\ref{lemrenewalestimates}}
The very first bound results simply, by monotone convergence, from $\Ep
[(M'_1)^\alpha(M_2)^\beta\ee^{\gamma H}]<\infty$ when $\gamma
<\kappa$,
which is a consequence, via H\"older inequality, of the fact that all
the moments of $M'_1$ and $M_2$ are finite under $\Pp$ (because
$M'_1\leq R_-+e_1$ and $M_2\leq e_1$), and of the fact that, due
to \eqref{iglehartthm}, $\ee^H$ has moments up to order $\kappa$ (not
included). Let us turn to the other bounds.

Note that, if $M'_1$ and $M_2$ were positive constants, then the bounds
would follow by an elementary computation from the tail estimate \eqref
{iglehartthm} and the classical formulas
\[
E\bigl[\ee^{\gamma H}\indic{H\geq h}\bigr] =\ee^{\gamma h}P(H\geq h)+
\int_h^\infty\gamma\ee^{\gamma u}P(H\geq u)\,
\d u
\]
and $E[\ee^{\gamma H}\indic{H<h}] =1-\ee^{\gamma h}P(H\geq h)+\int_0^h
\gamma\ee^{\gamma u}P(H\geq u)\,\d u$.

As recalled in Section~\ref{secsketch}, it was proved in \cite
{renewal} that indeed $M_1$ and $M_2$ depend little on $H$, in that
(Lemma 4.1 of~\cite{renewal}) for any integer\vadjust{\goodbreak} $r>0$ there is a constant
$C$ such that
%
%
\setcounter{equation}{0}
\begin{equation}
\label{eqnborneM1} \Ep\bigl[(M_1)^r |\lfloor H\rfloor,
H=S \bigr]\leq C,
\end{equation}
and similarly for $M_2$ (due to a symmetry property under $P^{\ge
0}(\cdot\vert H=S);$ see Lemma 3.4 in~\cite{renewal}). Admitting
that furthermore,
%
%
\begin{equation}
\label{eqnborneM12} \Ep\bigl[\bigl(M'_1
\bigr)^r |\lfloor H\rfloor, H=S \bigr]\leq C,
\end{equation}
we would first get by the Cauchy--Schwarz inequality that, with
$M\defeq\break
(M'_1)^\alpha(M_2)^\beta$,
\begin{eqnarray*}
&&\Ep\bigl[M|\lfloor H\rfloor, H=S\bigr] \\
&&\qquad\leq\Ep\bigl[\bigl(M'_1
\bigr)^{2\alpha}|\lfloor H\rfloor, H=S\bigr]^{1/2}\Ep
\bigl[(M_2)^{2\beta}|\lfloor H\rfloor, H=S\bigr]^{1/2}
\\
&&\qquad\leq C,
\end{eqnarray*}
and, using conditioning on $\lfloor H\rfloor$, conclude that
\[
\Ep\bigl[M\ee^{\gamma H}\indic{H<h}|H=S\bigr]\leq C' E\bigl[
\ee^{\gamma(\lfloor
H\rfloor+1)}\indic{\lfloor H\rfloor<h}\bigr]\leq C''
E\bigl[\ee^{\gamma
H}\indic{H<h+1}\bigr],
\]
and similarly
$
\Ep[M\ee^{\gamma H}\indic{H\geq h}|H=S]\leq C'' E[\ee^{\gamma
H}\indic
{H\geq h-1}]
$
which brings us back to the situation where $M'_1$ and $M_2$ would be constants.
Thus, it remains to prove \eqref{eqnborneM12} and, first, justify why
introducing the convenient condition $\{H=S\}$ is harmless.

As in Lemma~\ref{lemtailz}, the condition $\{H=S\}$ is typically
satisfied when $H$ is large; thus it suffices to note that the
contribution to the expectations of small values of~$H$ is not too
significant. Let $\ell=\ell(h)\defeq\frac{1}{\gamma}\log h$. We have
%
%
\begin{equation}
\label{eqn1464} \Ep\bigl[M \ee^{\gamma H} \indic{H<h}\bigr] \leq
\Ep[M] h+\Ep
\bigl[M \ee^{\gamma H} \indic{H<h, H>\ell}\bigr].
\end{equation}
Since $M$ and $H$ are independent of $S'\defeq\sup_{x\geq e_1}
V(x)-V(e_1)$, and also $\{S>H>\ell\}\subset\{S'>\ell\}$, we have on the
other hand
\[
\Ep\bigl[M \ee^{\gamma H}\indic{H<h, H>\ell}\indic{S>H}\bigr
]\leq\Ep\bigl[M
\ee^{\gamma H} \indic{H<h}\bigr]P\bigl(S'>\ell\bigr)
\]
and $P(S'>\ell)=o(1)$ when $h\to\infty$, hence, substracting this
quantity to \eqref{eqn1464} gives
\[
\Ep\bigl[M \ee^{\gamma H} \indic{H<h}\bigr]\bigl(1+o(1)\bigr
)\leq\Ep[M]h+\Ep
\bigl[M\ee^{\gamma
H} \indic{H<h}\indic{H=S}\bigr].
\]
Given that $P(H=S)>0$, and $h\leq\ee^{(\gamma-\kappa)h}$ for large $h$
when $\gamma>\kappa$, it thus suffices to prove the last two bounds of
\eqref{eqnm1m2inf} with the left-hand side replaced by $\Ep[M \ee
^{\gamma H}|H<h,H=S]$. As for \eqref{eqnm1m2sup}, the introduction of
$\ell$ is useless to similarly prove [skipping \eqref{eqn1464}] that
we may condition by $\{H=S\}$.

Let us finally prove \eqref{eqnborneM12}. Let $r>0$. We have
$M'_1=M_1+\sum_{T_H<x<e_1}\ee^{-V(x)}$. It results from Lemma 3.4 of
\cite{renewal} that $(H,\sum_{T_H\leq x<e_1}\ee^{-V(x)})$ has the same
distribution under $\Pp(\cdot|H=S)$ as $(H,\sum_{T^-_H<x\leq0}\ee
^{V(x)-H})$ where
\[
T^-_H\defeq\sup\bigl\{x\leq0 \dvtx V(x)>H\bigr\},
\]
and we claim that there is $C'_r>0$ such that, for all $N\in\N$,
%
%
\begin{equation}
\label{eqn4646} E \biggl[ \biggl(\sum_{T^-_N<x\leq0}
\ee^{V(x)} \biggr)^r \biggr]\leq C'_r
\ee^{rN}.
\end{equation}
Before we prove this inequality, let us use it to conclude that
%
%
\begin{eqnarray}\label{eqnbornem15}
&&\Ep\bigl[\bigl(M'_1\bigr)^r|\lfloor H
\rfloor, H=S\bigr]
\nonumber
\\[-8pt]
\\[-8pt]
\nonumber
&&\qquad\leq2^r \bigl(\Ep\bigl[(M_1)^r|
\lfloor H\rfloor, H=S\bigr]+\ee^{-r\lfloor
H\rfloor}C \ee^{r(\lfloor H\rfloor+1)} \bigr)\leq C'.
\end{eqnarray}
For readability reasons, we write the proof of \eqref{eqn4646} when
$r=2$, the case of higher integer values being exactly similar and
implying the general case (if $0<r<s$, $E[X^r]\leq E[X^s]^{r/s}$ for
any positive $X$).
We have
%
%
\begin{eqnarray}
\label{eqn45678}
&&E \biggl[ \biggl(\sum_{T_N^-<x\leq0}
\ee^{V(x)} \biggr)^2 \biggr]
\nonumber
\\[-8pt]
\\[-8pt]
\nonumber
&&\qquad\leq\sum
_{0\leq m,n<N} \ee^{n+1}\ee^{m+1}E\bigl[\nu\bigl([n,n+1
)\bigr) \nu\bigl([m,m+1)\bigr)\bigr],
\end{eqnarray}
where $\nu(A)\defeq\#\{x\leq0\dvtx V(x)\in A\}$ for all $A\subset
\R$. For
any $n\in\N$, applying the Markov property at time $\sup\{x\leq0\dvtx
V(x)\in
[n,n+1)\}$ gives us that
$E[\nu([n,n+1))^2]\leq E[\nu([-1,1))^2]$. This latter expectation is
finite because $V(1)$ has a negative mean and is exponentially
integrable; more precisely, $\nu([-1,1))$ is exponentially integrable
as well: for $\lambda>0$, for all $x\geq0$, $P(V(-x)<1)\leq\ee
^{\lambda}E[\ee^{\lambda V(1)}]^x=\ee^\lambda E[\rho^\lambda]^x$ hence,
choosing $\lambda>0$ small enough so that $E[\rho^\lambda]<1$ [cf.
Assumption (a)], we have, for all $p\geq0$,
\begin{eqnarray*}
P\bigl(\nu\bigl([-1,1\bigr)\bigr)>p) &\leq& P\bigl(\exists x\geq
p\mbox{ s.t. }
V(-x)<1\bigr)
\\
& \leq&\sum_{x\geq p}P\bigl(V(-x)<1\bigr)\leq
\ee^\lambda\bigl(1-E\bigl[\rho^\lambda\bigr]\bigr)^{-1}
E\bigl[\rho^\lambda\bigr]^p.
\end{eqnarray*}
Thus, using the Cauchy--Schwarz inequality to bound the expectations
uniformly, the right-hand side of \eqref{eqn45678} is less than $C
\ee^{2N}$ for some constant $C$. This proves~\eqref{eqn4646} and
therefore concludes the proof of Lemma~\ref{lemrenewalestimates}.
\end{pf*}

\subsection{\texorpdfstring{Proofs of Lemmas \protect\ref{lemgoodenv}, \protect\ref{lemmafailurebound} and \protect\ref{lemmafsuccessbound}}
{Proofs of Lemmas 5, 6 and 7}}
\label{subsecprooflemma-singlevalley}
\mbox{}
\begin{pf*}{Proof of Lemma~\ref{lemgoodenv}} By the union bound the proof
of Lemma~\ref{lemgoodenv} boils down to showing that for $i=1,2,3,$
\[
P\bigl(\bigl(\Omega_t^{(i)}\bigr)^c, H \ge
\hP_t\bigr)=o\bigl(t^{-\kappa}\bigr),\qquad t \to\infty.\vadjust{\goodbreak}
\]

The case $i=1$ is trivial. Indeed, the fact that $e_1$ has some finite
exponential moments (see Section~\ref{subsecrenewal}) implies that
$P((\Omega_t^{(1)})^c)=o(t^{-\kappa})$ when $t$ tends to infinity (for
$C$ large enough).
The case $i=2$ can be proved by a minor adaptation of the proof of
Lemma 5.5 in~\cite{limitlaws}.

Let us consider the last case $i=3$. Since $R^-$ depends only on the
variables $V(x), x \le0$, and $P(H>\hP_t)\sim C_I t^{-\kappa
}(\log
t)^{\kappa}$ when $t\to\infty$, it suffices to prove $\Pp(R^->(\log
t)^4 t^\alpha)=o((\log t)^{-\kappa})$. This would follow (for any
$\alpha>0$) from the Markov property if $\Ep[R^-]<\infty$. We have
(changing indices and incorporating the single terms into the sums)
%
%
\begin{eqnarray}
R^- &=& \sum_{x\leq0} \biggl(1+2\sum
_{x< y\leq0}\ee^{V(y)-V(x)} \biggr) \biggl(\ee^{-V(x)}+2
\sum_{z\leq x-1}\ee^{-V(z)} \biggr)
\nonumber
\\[-8pt]
\\[-8pt]
\nonumber
&\leq&4\sum_{z\leq x\leq y\leq0}\ee^{V(y)-V(x)-V(z)},\label{eqnr^-}
\end{eqnarray}
and this latter quantity was already seen to be integrable under $\Pp$,
after \eqref{eqnjneg}, when $1<\kappa<2$. In order to deal with the
case $0<\kappa\leq1$, let us introduce the event
\[
A_t=\bigcap_{p=1}^\infty
\biggl\{H_{-p}<\frac{1}{\kappa}\log p^2+\log t+\log\log
t\biggr\}.
\]
On one hand, by \eqref{iglehartthm}, $P((A_t)^c)\leq\sum
_{p=1}^\infty
\frac{C}{p^2 (t\log t)^\kappa}= (\sum_{p=1}^\infty\frac
{C}{p^2} )\frac{t^{-\kappa}}{(\log t)^\kappa}=o(t^{-\kappa
})$. On
the other hand, proceeding as after \eqref{eqnjneg},
\begin{eqnarray*}
&&\Ep\bigl[R^-\ind_{A_t}\bigr] \\
&&\qquad\leq4\sum_{p\leq0}
\Ep\bigl[\ee^{-V(e_p)}\bigr]\Ep\bigl[\bigl(M'_1
\bigr)^2 M_2 \ee^H \indic{H<
({1}/{\kappa})\log p^2+\log t+\log\log t}\bigr]
\end{eqnarray*}
and $\Ep[\ee^{-V(e_p)}]=E[\ee^{V(e_1)}]^p$ hence, using Lemma \ref
{lemrenewalestimates}, when $0<\kappa<1$,
\[
\Ep\bigl[R^-\ind_{A_t}\bigr] \leq4 \biggl(\sum
_{p\leq0}E\bigl[\ee^{V(e_1)}\bigr]^p
\frac
{1}{(p^2)^{(1-\kappa
)/\kappa}} \biggr) (t\log t)^{1-\kappa} \le C(t\log
t)^{1-\kappa},
\]
and when $\kappa=1$,
\[
\Ep\bigl[R^-\ind_{A_t}\bigr] \leq4\sum_{p\leq0}E
\bigl[\ee^{V(e_1)}\bigr]^p \biggl(\frac{1}{\kappa}\log
p^2+\log t+\log\log t \biggr)\leq C\log t.
\]
Finally, by the Markov inequality,
\[
\Pp\bigl(R^->t^\alpha(\log t)^4\bigr)\leq\Pp
\bigl((A_t)^c\bigr)+\frac{1}{t^\alpha
(\log
t)^4}\Ep\bigl[R^-
\ind_{A_t}\bigr]
\]
is negligible with respect to $(\log t)^{-\kappa}$ for any $\alpha
\geq
1-\kappa$ when $0<\kappa<1$, and for any $\alpha>0$ when $\kappa=1$.\vadjust{\goodbreak}
\end{pf*}

\begin{pf*}{Proof of Lemma~\ref{lemmafailurebound}}
Since $\Var_\omega(F) \le E_{\omega} [ F^2 ],$ the proof of
\eqref{eqboundVarF} is a consequence of (5.10) in~\cite{limitlaws}
together with a minor adaptation of equation (5.26) in~\cite{limitlaws}
and the definition of $\Omega_t.$
The proof of~\eqref{eqboundM2} is a direct consequence of the
definitions of $M_2$ [see equation \eqref{eqdefM1M2}] and $\Omega_t$
[see equation~\eqref{eqgoodenv}].\vspace*{1pt}
Finally, the proof of~\eqref{eqboundM1} is straightforward by looking
at the expression of $ E_{\omega} [ F ] = 2 \omega_0
\widehat M_1$ in terms of the modified potential $\widehat{V}$ (see
Lemma~5.2 in~\cite{limitlaws}) together with the properties of good
environments $\omega$ in~$\Omega_t.$
\end{pf*}

\begin{pf*}{Proof of Lemma~\ref{lemmafsuccessbound}} The proof of Lemma
\ref
{lemmafsuccessbound} can be deduced from Lemma~5.4 in~\cite{limitlaws}
(which gives an upper bound for $E_{\omega}[G]$ in terms of the
modified potential~$\bar V$), the definition of the modified potential
$\bar V$ (see equation (5.15) in~\cite{limitlaws}) and the definition
of good environments $\omega$ in $\Omega_t.$
\end{pf*}

\subsection{An annealed result}
The techniques of this paper enable us to prove the following annealed
counterpart to \eqref{eqqueueetau} which has its own interest.
%
\begin{proposition}\label{propplougonvelin}
The tail distribution of the hitting time of the first negative record
$e_1$ satisfies
%
%
\begin{equation}
t^\kappa\bPp\bigl(\tau(e_1)\geq t\bigr)\limites{} {}
C_T,\qquad t\to\infty,
\end{equation}
where the constant $C_T$ is given by
%
%
\begin{equation}
C_T\defeq2^\kappa\Gamma(\kappa+1)C_U.
\end{equation}
\end{proposition}

Let us write $\tau$ for $\tau(e_1)$ in this section. The idea of the
proof is the following. We first show that, on the event $\{\tau\geq
t\}
$, the height of the first excursion is typically larger than the
function $\hP_t$ [of order $\log t$, defined in~\eqref{hTdef}]. We
may then invoke Proposition~\ref{propexitonevalley} to reduce the
tail of $\tau$ to that of $E_\omega[\tau]\eg$ and conclude.

%
\begin{lemma} \label{lemHlarge} We have
\[
\label{eqestimatevalley1} \bPp\bigl(\tau(e_1)\geq t , H <
\hP_t\bigr)=o\bigl(t^{-\kappa}\bigr),\qquad t \to\infty.
\]
\end{lemma}

\begin{pf}
Let us first assume that $0<\kappa<1$. Then, by the Markov inequality, we get
\begin{eqnarray*}
\bPp(\tau\geq t, H< \hP_t) & =& \Ep\bigl[P_\omega(\tau\geq
t)\indic{H< \hP_t}\bigr] \leq\frac{1}{t} \Ep
\bigl[E_\omega[\tau]\indic{H< \hP_t}\bigr]
\\
& \leq&\frac{1}{t} \Ep\bigl[2M_1'M_2
\ee^H \indic{H< \hP_t}\bigr] \leq\frac
{1}{t} C
\ee^{(1-\kappa)\hP_t},
\end{eqnarray*}
where the last inequality follows from Lemma \ref
{lemrenewalestimates}. Since we have $t^{-1}\ee^{(1-\kappa)\hP
_t}=t^{-\kappa}(\log t)^{-(1-\kappa)}$, this settles this case.

Let us now assume $1<\kappa<2$. By the Markov inequality, we get
\[
\bPp(\tau\geq t, H <\hP_t) \le\frac{1}{t^2} \Ep
\bigl[E_\omega\bigl[\tau^2\bigr]\indic{H <
\hP_t}\bigr].
\]
Applying Lemma~\ref{lemtaue1} yields $\bPp(\tau\geq t, H <\hP_t)
\le Ct^{-2} \ee^{(2-\kappa)\hP_t}$, which concludes the proof of Lemma
\ref{lemHlarge} when $\kappa\neq1$.

For $\kappa=1$, neither of the above techniques works; the first one is
too rough, and $\Var_\omega(\tau)$ is not integrable hence, the second
does not make sense as is. We shall modify $\tau$ so as to make $\Var
_\omega(\tau)$ integrable. To this end, let us refer to
Section~\ref
{subsecreducsmall} and denote by $d_-$ the right end of the first
excursion on the left of 0 that is higher than $h_t$, and by $\taut
\defeq\tautx{d_-}(0,e_1)$ the time spent on the left of $d_-$ before
reaching~$e_1$. By Lemma~\ref{lemremonteedifficile} we have $\bEp
[\taut\indic{H<\hP_t}]\leq C \hP_t\ee^{-\hP_t}\leq C(\log t)^2t^{-1}$.
Let us also introduce $\taut'$, which is defined like $\taut$ but in
the modified environment, that is, by replacing the high excursions
(on the left of $d_-$) by small ones; cf. after Lemma~\ref
{lemremonteedifficile}. Then we have
\begin{eqnarray*}
& &\bPp(\tau\geq t,H<\hP_t)
\\
& &\qquad\leq\bPp\bigl(\taut\geq(\log t)^3,H<\hP_t\bigr)+\bPp
\bigl(\tau-\taut\geq t-(\log t)^3, H<\hP_t\bigr)
\\
&&\qquad\leq\frac{1}{(\log t)^3}\bEp[\taut\indic{H<\hP_t}] +
\bPp\bigl(
\tau-\taut+\taut' \geq t-(\log t)^3, H<
\hP_t\bigr)
\\
&&\qquad= o\bigl(t^{-1}\bigr) + \bigl(\bPp\bigr)'\bigl(\tau\geq t-(\log
t)^3, H<\hP_t\bigr)
\\
&&\qquad\leq o\bigl(t^{-1}\bigr) + \frac{1}{(t-(\log t)^3)^2}\bigl(\Ep\bigr)'
\bigl[E_\omega\bigl[\tau^2\bigr]\indic{H<\hP_t}
\bigr],
\end{eqnarray*}
and Lemma~\ref{lemtaue1} allows us to conclude just like in the case
$1<\kappa<2$.
\end{pf}

\begin{pf*}{Proof of Proposition~\ref{propplougonvelin}}
From the tail of $E_\omega[\tau]$ [cf. \eqref{eqqueueetau}], a simple
computation gives
%
%
\begin{equation}
\label{eqntailEe} \Pp\bigl(E_\omega[\tau]\eg\geq t \bigr)\sim
C_T t^{-\kappa},\qquad t\to\infty.
\end{equation}
Let us prove that this is also the tail of $\tau$.

For any function $t\mapsto u_t$ we have, using, respectively, the
previous lemma for the first bound and Proposition \ref
{propexitonevalley} and the Markov inequality (with respect to $P_\omega
$) for the second,
\begin{eqnarray*}
\bPp\bigl(\tau-E_\omega[\tau]\eg\geq u_t, \tau>t\bigr) &
\leq&\bPp\bigl(\tau-E_\omega[\tau]\eg\geq u_t, H\geq
\hP_t\bigr)+o\bigl(t^{-\kappa
}\bigr)
\\
& \leq&\frac{t^\beta}{u_t}\bP(H\geq\hP_t)+o\bigl(t^{-\kappa}
\bigr)
\\
&=&t^{-\kappa} \biggl(\frac{t^\beta\log t}{u_t}\bigl(1+o(1)\bigr
)+o(1) \biggr).
\end{eqnarray*}
If we choose $u_t$ such that $t^\beta(\log t)^\kappa\ll u_t\ll t$ then
we get, assembling this with~\eqref{eqntailEe},
\begin{eqnarray*}
\bPp(\tau>t) & =& \bPp\bigl(\tau-E_\omega[\tau]\eg\geq
u_t,\tau>t\bigr)+\bPp\bigl(\tau-E_\omega[\tau]\eg<
u_t,\tau>t\bigr)
\\
& \leq& o\bigl(t^{-\kappa}\bigr)+\bPp\bigl(E_\omega[\tau]\eg
\geq
t-u_t\bigr)\sim C_Tt^{-\kappa}.
\end{eqnarray*}
The lower bound is identical, starting with
\[
\bPp(\tau> t)\geq\bPp\bigl(E_\omega[\tau]\eg\geq t+u_t
\bigr)-\bPp\bigl(\tau-E_\omega[\tau]\eg\leq-u_t,\tau>t
\bigr) .
\]
This concludes the proof of Proposition~\ref{propplougonvelin}.
\end{pf*}
\end{appendix}


\section*{Acknowledgment}
Many thanks are due to an anonymous referee for
careful reading of the original manuscript and helpful comments.

%

%

%


\printaddresses

\end{document}